\documentclass{amsart}
\usepackage[left=3.3cm,right=3.2cm]{geometry}
\usepackage{bbm,pifont}
\usepackage[hidelinks]{hyperref}
\usepackage{amsmath,amssymb,amscd,amsthm,latexsym,mathrsfs}
\usepackage{graphicx,xcolor,subfigure,extarrows,caption}
\usepackage{pdfpages}
\usepackage{url}
\usepackage[english]{babel}
\usepackage{mathtools}
\usepackage{comment}

\newtheorem{thm}{Theorem}[section]
\newtheorem{lemma}[thm]{Lemma}
\newtheorem{cor}[thm]{Corollary}
\newtheorem{prop}[thm]{Proposition}

\theoremstyle{definition}
\newtheorem{defi}[thm]{Definition}

\newtheorem*{rmk}{Remark}

\newcommand{\gauss}[1]{{\lfloor{#1}\rfloor}}

\newcommand{\EC}{\widehat{\mathbb{C}}}
\newcommand{\C}{\mathbb{C}}

\newcommand{\BH}{\mathbb{H}}
\newcommand{\BL}{\mathbb{L}}

\newcommand{\J}{\mathbb{J}}
\newcommand{\N}{\mathbb{N}}
\newcommand{\Q}{\mathbb{Q}}
\newcommand{\R}{\mathbb{R}}

\newcommand{\Z}{\mathbb{Z}}

\newcommand{\MA}{\mathcal{A}}

\newcommand{\MC}{\mathcal{C}}
\newcommand{\MD}{\mathcal{D}}
\newcommand{\ME}{\mathcal{E}}
\newcommand{\MF}{\mathcal{F}}
\newcommand{\MG}{\mathcal{G}}
\newcommand{\MH}{\mathcal{H}}

\newcommand{\MK}{\mathcal{K}}
\newcommand{\ML}{\mathcal{L}}

\newcommand{\MO}{\mathcal{O}}
\newcommand{\MP}{\mathcal{P}}
\newcommand{\MQ}{\mathcal{Q}}
\newcommand{\MMR}{\mathcal{R}}

\newcommand{\IS}{\mathcal{IS}}
\newcommand{\QIS}{\mathcal{QIS}}

\newcommand{\HB}{\mathscr{B}}

\newcommand{\HH}{\mathscr{H}}
\newcommand{\HJ}{\mathscr{J}}
\newcommand{\HS}{\mathscr{S}}

\newcommand{\re}{\operatorname{Re}}
\newcommand{\im}{\operatorname{Im}}

\newcommand{\Expo}{\operatorname{\mathbb{E}\textup{xp}}}

\newcommand{\diam}{\operatorname{diam}}

\newcommand{\Boxx}{\operatorname{Box}}
\newcommand{\cv}{\textup{cv}}
\newcommand{\cp}{\textup{cp}}
\newcommand{\HT}{\textup{HT}}

\newcommand{\Dom}{\operatorname{Dom}}
\newcommand{\Int}{\operatorname{int}}

\newcommand{\dens}{\operatorname{dens}}
\newcommand{\area}{\operatorname{area}}
\newcommand{\Pack}{\operatorname{Pack}}

\newcommand{\kc}{\textit{\textbf{k}}}
\newcommand{\ea}{\textup{\textbf{s}}}

\newcommand{\abs}[1]{\left| #1 \right|}
\newcommand{\p}[1]{\left( #1 \right)}

\makeatletter\@addtoreset{equation}{section}\makeatother

\begin{document}

%---------------------------------------------------------------------------------------------------------------
\title[Dimension paradox of irrationally indifferent attractors]{Dimension paradox of irrationally indifferent attractors}

\author[D. Cheraghi]{Davoud Cheraghi}
\address{Department of Mathematics, Imperial College London, London SW7 2AZ, UK}
\email{d.cheraghi@imperial.ac.uk}

\author[A. DeZotti]{Alexandre DeZotti}
\address{Department of Mathematics, Imperial College London, London SW7 2AZ, UK}
\email{a.de-zotti@imperial.ac.uk}

\author[F. Yang]{Fei YANG}
\address{School of Mathematics, Nanjing University, Nanjing 210093, P. R. China}
\email{yangfei@nju.edu.cn}

\begin{abstract}
%In this paper we study the geometry of the post-critical set for a class of holomorphic maps with a neutral fixed point. 
We prove that for an infinite dimensional class of holomorphic maps with an elliptic fixed point, 
the post-critical set has Hausdorff dimension two, provided the rotation number is a non-Herman of sufficiently high type.
%As an immediate corollary, we conclude that the Hausdorff dimension of the Julia set of a class of rational functions with 
%a Cremer fixed point is equal to two.
We identify classes of Brjuno but not Herman, and non-Brjuno, numbers such that the post-critical set satisfies 
the Karpi\'{n}ska's dimension paradox. 
That is, the set of end points of the post-critical set has dimension two, but without those end points, the dimension drops to one.
\end{abstract}

% AMS subject classifications (used in AMS journals)

% AMS keywords (used in AMS journals)

% today's date, or fill in whatever date you prefer
\date{\today}

% acknowledge support, etc
% \thanks{This research was partially supported by NSF grant DOA-123456789.}

% dedication
% \dedicatory{Dedicated to Professor}

\maketitle

%----------------------------------------------------------------------------------------------------------------

\section{Introduction}\label{introduction}
Consider the quadratic polynomials 
\[P_\alpha(z)=e^{2\pi i \alpha}z+z^2\]
on the complex plane $\C$, where $\alpha$ is an irrational number with continued fraction 
\[\alpha= 1/(a_{0}+1/(a_1+1/(a_2+\dots ))).\] 
The \textit{post-critical} set of $P_\alpha$ is the closure of the orbit of the unique critical point of $P_\alpha$, % $-e^{2\pi i \alpha}/2$,   
\[\Lambda(P_\alpha)= \overline{\textstyle{\bigcup_{i\geq 1}} P_\alpha^{\circ i}(-e^{2\pi i \alpha}/2)}.\]

For rotation numbers of \textit{bounded type}, that is, when $\sup_{i\geq 0} a_i < \infty$, 
$\Lambda(P_\alpha)$ is well understood. 
By the work of Siegel \cite{Sie42}, $P_\alpha$ is conformally conjugate to rotation by $2\pi \alpha$ near $0$. 
Using a surgery procedure of Douady~\cite{Do86}, %sing the maximal linearisation domain, an ingenious surgery procedure by Douady \cite{Do86} 
%links the dynamics of $P_\alpha$ to the dynamics of certain cubic Blasche products, where the seminal works of 
Herman and Swianek \cite{He87,Swi98} showed that $\Lambda(P_\alpha)$ is a \textit{quasi-circle}. 
Thus, $\Lambda(P_\alpha)$ encloses the maximal invariant open set on which $P_\alpha$ is conjugate to the rotation, 
now called the \textit{Siegel disk} of $P_\alpha$.
Roughly speaking, a quasi-circle is a Jordan curve with controlled geometry at all scales. 
By a renormalisation method, McMullen \cite{McM98} showed that $\Lambda(P_\alpha)$ enjoys a rescaling self-similarity 
at the critical point when $(a_i)_{i\geq 0}$ is eventually periodic. 
Generalising the surgery procedure, Petersen and Zakeri \cite{PZ04} showed that under the growth condition $\log a_i = \mathcal{O}(\sqrt{i})$, 
$\Lambda(P_\alpha)$ is a \textit{David circle} (a Jordan curve with controlled degenerating geometry at small scales). 
%It is known that large entries in the continued fraction stretches the geometry of $\Lambda(P_\alpha)$. 
%However, the relation has not been fully understood.  

When an entry $a_i$ in the continued fraction of $\alpha$ tends to infinity, $\alpha$ tends to a rational number, 
and the dynamic of $P_\alpha$ degenerates to a parabolic one. 
When some or all entries are simultaneously very large, one has a plethora of ``near degenerations''. 

Using the near-parabolic renormalisation scheme of Inou and Shishikura \cite{IS06}, 
a complete topological description of $\Lambda(P_\alpha)$ is recently established in \cite{Che17} 
when $\alpha$ is of sufficiently \textit{high type}. 
That is, for irrational numbers $\alpha$ with $\inf_{i\geq 0} a_i \geq N$, for some specific $N$ 
determined by the renormalisation scheme. 
We denote those numbers by $\HT_N$.  
For every $\alpha \in \HT_N$, one of the following holds: 
\begin{itemize}
\item[(i)] $\alpha$ is a Herman number, and $\Lambda(P_\alpha)$ is a Jordan curve, 
\item[(ii)] $\alpha$ is a Brjuno but not Herman number, and $\Lambda(P_\alpha)$ is a \textit{one-sided hairy Jordan curve},
\item[(iii)] $\alpha$ is not a Brjuno number, and $\Lambda(P_\alpha)$ is a \textit{Cantor bouquet}.
\end{itemize}

%The topological classification is determined by the arithmetic type of Herman and the more relaxed 
%arithmetic type of Brjuno. 
For the definitions of the Herman and Brjuno types see Section~\ref{section:arithmetic rot numbers}, 
and for the definitions of one-sided hairy Jordan curve and Cantor bouquet, see Section \ref{subsec:topo}. 
Any Herman number is a Brjuno number.
Roughly speaking, in case (iii), $\Lambda(P_\alpha)$ is a collection of arcs landing at $0$ where 
each arc in $\Lambda(P_\alpha)$ is accumulated from both sides by arcs in $\Lambda(P_\alpha)$. 
In case (ii), $\Lambda(P_\alpha)$ is a collection of arcs landing on one side of a Jordan curve where 
each arc in the collection is accumulated from both sides by arcs in $\Lambda(P_\alpha)$.
%In cases (i) and (ii), the Jordan curve in $\Lambda(P_\alpha)$ is the boundary of the Siegel disk of $P_\alpha$. 
See also \cite{FSh18} for a partial result towards the above trichotomy.  \nocite{GrSw03}

As in the bounded type regime, $\Lambda(P_\alpha)$ is the frontier separating the quasi-periodic behaviour 
from the global chaotic behaviour.
More precisely, the filled-in $\Lambda(P_\alpha)$ is the maximal connected invariant set containing $0$ 
where $P_\alpha$ acts quasi-periodically. 
In particular, in cases (i) and (ii) above, the region enclosed by the Jordan curve in $\Lambda(P_\alpha)$
is the Siegel disk of $P_\alpha$.

Despite the trichotomy of topology, a universe of geometric possibilities is present. 
Here we study the geometry of $\Lambda(P_\alpha)$ when degenerations have occurred in cases (ii) and (iii). 

\begin{thm}\label{thm:quadratic-dim-pc}
For every non-Herman $\alpha$ in $\HT_N$, $\Lambda(P_\alpha)$ has Hausdorff dimension $2$.
\end{thm}

In \cite{Che13,Che19}, it is proved that for every $\alpha \in \HT_N$, $\Lambda(P_\alpha)$ has zero area 
but is the limit set of the orbit of almost every point in the Julia set of $P_\alpha$.  
On the other hand, Buff and Cheritat in \cite{BC12} proved the existence of $\alpha \in \HT_N$ 
such that the Julia set of $P_\alpha$ has positive area.
It remains open for which values of $\alpha$ the Julia set of $P_\alpha$ has positive area. 

\begin{cor}\label{cor:quadratic-dim-2}
For all non-Herman $\alpha$ in $\HT_N$, the Julia set of $P_\alpha$ has Hausdorff dimension $2$. 
\end{cor}

Shishikura in \cite{Sh98} proved that for a residual set of $\alpha$ in $\R/\Z$ the Julia set of $P_\alpha$ has 
Hausdorff dimension 2. But an arithmetic characterisation leading to that result was not available. 
For bounded type $\alpha$, the Hausdorff dimension of the Julia set of $P_\alpha$ is strictly less than $2$, \cite{McM98}.  

In this paper we present examples of $\Lambda(P_\alpha)$ with a peculiar geometric feature. 
%known as the dimension paradox of Karpi\'{n}ska. 
Let $\ME(P_\alpha)$ denote the set of all the end points of $\Lambda(P_\alpha)$, and let 
$\Delta(P_\alpha)$ denote the Siegel disk of $P_\alpha$, when it exists. 
%Recall that $\partial \Delta(P_\alpha)$ is the region inside the unique Jordan curve in $\Lambda(P_\alpha)$ in cases (i) and (ii) above. 
Then, in case (ii), $\Lambda(P_\alpha) \setminus \partial \Delta(P_\alpha)$ consists of uncountably many Jordan arcs (hairs), 
and in case (iii), $\Lambda(P_\alpha) \setminus \{0\}$ consists of uncountably many Jordan arcs (hairs). 

\begin{thm}\label{thm:dim-hairs-quadratic-(ii)}
There is a class of non-Herman but Brjuno irrational numbers $\HS$ such that for every $\alpha \in \HS$, 
\[\dim_H  (\Lambda(P_\alpha) \setminus (\partial \Delta(P_\alpha) \cup \ME(P_\alpha)))=1, \qquad \text{and} \qquad  
\dim_H \ME(P_\alpha)=2.\]
\end{thm}

\begin{thm}\label{thm:dim-hairs-quadratic-(iii)}
There is a class of non-Brjuno irrational numbers $\HJ$ such that for every $\alpha \in \HJ$,  
\[\dim_H (\Lambda(P_\alpha) \setminus (\ME(P_\alpha)))=1, \qquad \text{and} \qquad \dim_H \ME(P_\alpha)=2.\] 
\end{thm}

\begin{comment}
\begin{thm}\label{thm:dim-hairs-quadratic}
There are classes of irrational numbers $\HJ$ and $\HS$, such that elements of $\HJ$ are Brjuno but not Herman and 
elements of $\HS$ are non-Brjuno, satisfying the following:
\begin{itemize}
\item for every $\alpha \in \HJ$, $\dim_H  (\Lambda(P_\alpha) \setminus (\partial \Delta(P_\alpha) \cup \ME(P_\alpha)))=1$ and 
$\dim_H(\ME(P_\alpha))=2$; 
\item for every $\alpha \in \HS$, $\dim_H (\Lambda(P_\alpha) \setminus (\ME(P_\alpha)))=1$ and $\dim_H(\ME(P_\alpha))=2$. 
\end{itemize}
\end{thm}
\end{comment}

Theorems~\ref{thm:dim-hairs-quadratic-(ii)} and \ref{thm:dim-hairs-quadratic-(iii)} are surprising; 
the set of end points of a collection of disjoint arcs occupies more space 
than the set of those arcs without their end points. 
This remarkable paradox was first observed by Karpi\'{n}ska in the dynamics of the 
exponential maps $\lambda e^z$, for $0<\lambda< 1/e$, \cite{Kar99a,Kar99b}. 
Both classes $\HJ$ and $\HS$ are uncountable, and are determined by explicit arithmetic conditions 
in Section~\ref{section:arithmetic rot numbers}. 
%For example, if $a_i \to \infty$ sufficient fast, $\alpha \in \HJ$, while elements of $\HS$ are more difficult to come by. 

%In those papers, the especial form of the exponential map plays a prominent role, while in this paper, 
%we exploit the relations between the arithmetic of $\alpha$ and the nonlinearities of the large iterates of $P_\alpha$.

For the proofs, we make use of the near-parabolic renormalisation of Inou and Shishikura \cite{IS06}. 
In particular, the above theorems apply to an infinite dimensional class of maps of the form $e^{2\pi i \alpha}z+\MO(z^2)$, 
such as some cubic polynomials. 
See Theorems \ref{thm:dim-pc} and \ref{thm:dim-hairs}, as well as Corollary~\ref{cor:dim-2}, for the more general statements. 
Also, strictly speaking, we only require the entries in a modified notion of continued fractions to be of sufficiently 
high type. That allows infinitely many small entries in the standard continued fraction of $\alpha$. 
Despite that, irrationally indifferent attractors remain mysterious for arbitrary polynomials and rational functions 
in the unbounded regime.
In contrast, for all polynomials and rational functions, bounded type Siegel disks are quasi-disks, \cite{Za99,Zha11}.

%If the boundary of the Jordan domain $U$ is $C^1$-smooth, P\'{e}rez-Marco has shown that such a compact set $K$ is
%in fact unique \cite{Per96}. He proved also that the non-linearizable hedgehogs (i.e., $0$ is a Cremer point) have no
%interior \cite{Per96}. The structure of such hedgehogs is topologically complex. For example, P\'{e}rez-Marco showed
%in \cite{Per94} that they are not locally connected at any point different from the fixed point.

\tableofcontents 

\vspace*{-1.2cm}

\subsection*{Frequently used notations} \hfill
\begin{itemize}
\item $\N$, $\Z$, $\Q$, $\R$, $\C$ and $\EC$ denote the set of all natural numbers (including $0$), integers, rational numbers, 
real numbers, complex numbers, and the Riemann sphere, respectively, 
\item given $a\in \C$, $Z \subseteq \C$ and $r>0$, $B_r(a)=\{z\in \C : |z-a|< r\}$, and $B_r(Z):=\cup_{z \in Z}B_r(z)$, 
\item for $a \in \C$ and $Z, W \subseteq \C$, $aZ=\{az:z\in Z\}$, $Z + a=\{z +a : z \in Z\}$, 
and $Z+W=\{z+w:z\in Z, w\in W\}$, 
\item for $Z \subseteq \C$, $\diam(Z) = \sup_{a,b \in Z} |a-b|$, 
\item $\Dom (f)$ denotes the domain of definition of a map $f$, 
\item for $y \in \R$, $\mathbb{L}_y=\{z \in \C : \im z = y \}$ and $\mathbb{H}_y=\{z \in \C : \im z \geq y\}$. 
\end{itemize}
%For $x\geq 0$, $\alphauss{x}$ denotes the integer part of $x$. 
\section{Brjuno, Herman, Jagged and Spiky numbers} \label{section:arithmetic rot numbers}
We work with a slightly modified notion of continued fractions that is more suitable for employing 
renormalisation methods. 
The modified continued fraction algorithm is defined as follows.
For $x \in \mathbb{R}$, let $d(x, \mathbb{ Z })= \min \{ | x - n |, n \in \mathbb{ Z } \} \in[0,1/2]$.
Fix an irrational number $\alpha \in (-1/2, 1/2)$, and let 
\[\alpha_0 = d (\alpha, \mathbb{ Z } ).\]
There is a unique $ \varepsilon_0 \in \{ \pm 1 \} $ such that
$ \alpha = \varepsilon_0 \alpha_0 $.
We define the sequence $(\alpha_n)_{ n \geq 0 }$ according to
\[\alpha_{ n + 1 } = d (1 / \alpha_n, \mathbb{ Z } ),\]
and then identify $a_n \in \mathbb{Z}$ and $\varepsilon_{n+1} \in \{\pm 1 \}$ such that
\[1/\alpha_n = a_n + \varepsilon_{ n + 1 } \alpha_{ n + 1 }.\]
It follows that $0< \alpha_n<1/2$ and $a_n \geq 2 $, for all $n \geq 0$.
Throughout we shall use the notation 
\[\alpha= [(a_0, \varepsilon_0),(a_1, \varepsilon_1),(a_2, \varepsilon_2), \dots]= 
\varepsilon_0/(a_0+ \varepsilon_1/(a_1+ \varepsilon_2/(a_2 + \dots))).\]
%These sequences provide the continued fraction in Equation \eqref{equ:high-type}.

Let $ \beta_{ -1 } = 1$ and $\beta_n = \prod_{i=0}^n \alpha_i$, for $ n \geq 0 $.
Yoccoz in \cite{Yoc95} introduced the Brjuno function 
\[ \mathcal{ B } (\alpha) = \textstyle{\sum_{ n = 0 }^\infty} \beta_{ n -1 } \log (1/\alpha_n).\]
%This is defined for irrational values of $\alpha$. 
An irrational number $\alpha$ is called a Brjuno number if $\mathcal{ B } (\alpha) < \infty$. 
%\end{equation}
%He also showed that the difference 
%\begin{equation}
%$\left| \mathcal{ B } (\alpha) - \sum_{ n = 1 }^{ \infty } (1/q_n) \log q_{ n + 1 } \right|$
%\end{equation}
%is uniformly bounded, independent of $\alpha$.
By the seminal works of Siegel-Brjuno-Yoccoz \cite{Sie42,Brj71,Yoc95}, $P_\alpha$ is linearisable 
iff $\alpha$ is a Brjuno number. 
In the setting of high-type irrational numbers, this optimality also  follows from the trichotomy for
$\Lambda(P_\alpha)$, Theorem~\ref{thm-Lambda-topo}. 

In \cite{Her79,Yoc02}, Herman and Yoccoz identified the optimal arithmetic condition for the linearisation of 
orientation preserving analytic circle diffeomorphisms. 
However, the arithmetic condition is only presented in terms of the standard continued fraction algorithm.
Below we present this arithmetic condition in terms of the modified continued fraction algorithm. 
The equivalence of the two conditions is proved in \cite{Che17}. 

For $r\in (0, 1/2)$, define the function $ h_r: \R \to \R$ as
\[h_r (y) =
\begin{cases}
r^{-1} (y+1- \log r^{-1})   &      \text{ if } y \geq \log r^{-1}, \\
e^y                                     &      \text{ if }  y \leq \log r^{-1}.
\end{cases}\]
The function $ h_r $ is $ C^1 $ and satisfies
\[\begin{gathered}
h_r (\log r^{-1}) =  h_r' (\log r^{-1}) = r^{-1}; \\
e^y \geq h_r (y) \geq y + 1, \forall y \in \mathbb{ R }; \\
h_r' (y) \geq 1, \;\;  \forall y \geq 0.
\end{gathered}
\]

An irrational number $\alpha$ is of \textbf{Herman type}, if for any $ n \geq 0 $ there is $ p \geq 1 $ such that
\[h_{ \alpha_{ n + p-1 } } \circ \cdots \circ h_{ \alpha_n } (0) \geq \mathcal{ B } (\alpha_{ n + p }).\]

In particular, any irrational number of Herman type belongs to $\mathscr{B}$.

Below, we define two classes of irrational numbers for which the conclusions of 
Theorems~\ref{thm:dim-hairs-quadratic-(ii)} and \ref{thm:dim-hairs-quadratic-(iii)} hold. 
For $ x \geq 0$, let
\[\lfloor x \rfloor = \max \{ n \in \mathbb{ \N } :  n \leq x \}\]
denote the integer part of $x$. 

An irrational number $\alpha= [(a_0, \varepsilon_0),(a_1, \varepsilon_1),(a_2, \varepsilon_2), \dots]$ 
is called a \textbf{jagged} number, if $\varepsilon_n=-1$ for all $n\geq 0$, and there is a sequence of positive numbers 
$(u_n)_{n\geq 0}$ such that
\begin{itemize}
\item[(i)] $\sum_{ n \geq 0 } u_0 \cdots u_n = + \infty$;
\item[(ii)] for all $n\geq 0$, $a_{n+1} \geq a_n^{u_n a_n}+1/2$;
\item[(iii)] $\lim_{n\to \infty} a_n = +\infty$; and
\item[(iv)] $u_n \log a_n \to + \infty$ as  $n\to\infty$.
\end{itemize}

The set of irrational numbers in $\HT_N$ which satisfies the above definition is denoted by $\HJ$.
For example, if $a_0=N$ and $a_{n+1}=\lfloor e^{e^{a_n}}\rfloor$, then $\alpha$ is a jagged number. 

\begin{lemma}
Any jagged number is a non-Brjuno number. 
\end{lemma}

\begin{proof}
By the choice of $\varepsilon_n$, for all $ n \geq 0 $, $ a_n - 1/2 < 1/ \alpha_n < a_n$.
Therefore, 
\[1/\alpha_{ n + 1 } > a_n^{ u_n a_n } > ( 1/ \alpha_n )^{u_n/\alpha_n}.\]
Then, for all $ n \geq 0 $,
$ \log (1/\alpha_{ n + 1 }) \geq (u_n/\alpha_n) \log (1/\alpha_n) $, and hence 
\[ \alpha_0 \cdots \alpha_n \log ( 1/ \alpha_{ n + 1 }) \geq u_n \alpha_0 \cdots \alpha_{ n - 1 } \log ( 1 / \alpha_{ n } ).\]
Then, $\sum_{ n \geq 0 } \beta_{ n - 1 } \log ( 1 / \alpha_n ) \geq 
(1+\sum_{ n \geq 0 } u_{ 0 } \cdots u_n)  \log ( 1 / \alpha_0) = + \infty$. 
\end{proof}

An irrational number $\alpha= [(a_0, \varepsilon_0),(a_1, \varepsilon_1),(a_2, \varepsilon_2), \dots]$ 
is called a \textbf{spiky} number if $\varepsilon_n=-1$ for all $n\geq 0$ and there are a sequence of 
positive numbers $(v_n)_{n \geq 0}$ and a uniformly bounded sequence of real numbers 
$(\eta_n)_{n \geq 0} $ such that
\begin{itemize}
\item[(i)] $v_n \to + \infty$, as $n \to + \infty$;
\item[(ii)] for all $n\geq 0$, $a_{ n + 1 } = e^{ v_n a_n } + \eta_n$; and
\item[(iii)] $\sum_{n \geq 1} v_n/(a_0 \cdots a_{ n - 1 }) < +\infty$.
\end{itemize}

The set of irrational numbers in $\HT_N$ which satisfies the above definition is denoted by $\HS$. 
For example, if $ (a_n)_{n \geq 1}$ satisfies $a_0=N$ and $a_{n+1} = \lfloor e^{2^n a_n} \rfloor+1$, then $\alpha$ is a spiky number.

\begin{lemma} \label{lem:spiky}
Any spiky number is of Brjuno type, but not of Herman type.
\end{lemma}

\begin{proof}
By the choice of $\varepsilon_n$, for all $ n \geq 0 $, we have 
\[ e^{ v_n / \alpha_n } + \eta_n - 1/2 \leq 1 / \alpha_{ n + 1 }  \leq e^{ v_n (1 / \alpha_n + 1 / 2) } + \eta_n.\]
Hence
\[ \alpha_0 \cdots \alpha_n \log  (1 / \alpha_{ n + 1 })  \leq \beta_{ n - 1 } \big(v_n (1 + \alpha_n / 2)  + \log (1 + C_n) \big),\]
with $ C_n \to 0 $ as $ n \to \infty $.
Then, there exists a constant $M>0$ such that
\[\mathcal{B} (\alpha) 
< \log (1/ \alpha_0) + (3/2) v_0 + (3/2) \textstyle{\sum_{ n \geq 1 }} ( v_n / a_0 \cdots a_{ n - 1 } ) + M <+\infty.
\]
On the other hand, since $v_n\to+\infty$ as $n\to\infty$, there exists 
$ n_0 \geq 0 $ such that for all $ n \geq n_0 $, $1 / \alpha_{ n + 1 } \geq e^{ 2 / \alpha_n }$.
In order to show that $ \alpha \notin \mathscr{ H } $ it is sufficient to show that for all $ n \geq n_0 $ 
and all $ p \geq 0 $, $ E^{\circ p} (0) < \log ( 1 / \alpha_{ n + p } )$, where $ E^{\circ p} $ is the $ p $-th 
iterate of the exponential map $ x \mapsto e^x $. 

Note that for $ p \geq 1 $, $\log (1 / \alpha_{ n + p } ) \geq  2 / \alpha_{ n + p - 1 } $.
In particular, we have 
\[ E (0) < 3 \leq 2 / \alpha_n \leq \log ( 1 / \alpha_{ n + 1 } ).\]
Moreover
\[ E^{\circ 2} (0) < 2 e^{ 2 / \alpha_n } \leq 2 /\alpha_{ n + 1 }  \leq \log (1 / \alpha_{ n + 2 }).\]
Similarly, using the function $E_2(x)= 2 e^x$, one can inductively prove that for $p \geq 1$, 
\[ E^{\circ p} (0) < E_2^{\circ (p - 1)} (2 / \alpha_n) \leq 2 / \alpha_{ n + p - 1 } \leq \log ( 1 / \alpha_{ n + p } ). \qedhere\]
\end{proof}

The terms, jagged and spiky, reflects the geometric features of the renormalisation towers associated with
such rotation numbers. This will be discussed in Section~\ref{section:HD-hairs}.

\begin{rmk}
It is possible to include rotation numbers with $\varepsilon_n=+1$ in the definitions of the jagged and spiky rotation numbers. 
Indeed, if $\alpha$ is a jagged or spiky rotation number with entries $(a_n, \varepsilon_n)_{n\geq 0}$, 
one may arbitrarily change the signs $\varepsilon_n$ to $\pm1$, while keeping the entries $a_n$ unchanged, and thereby obtain 
new irrational numbers which satisfy the conclusion of Theorems~\ref{thm:dim-hairs-quadratic-(ii)} and \ref{thm:dim-hairs-quadratic-(iii)}.  
However, since any positive value of $\varepsilon_n$ gives rise to an anti-holomorphic change of coordinates in the renormalisation tower, 
we have excluded them, for the sake of simplicity. 
\end{rmk}
\section{A criterion for full Hausdorff dimension}\label{sec-criterion}
In this section we present a general criterion by McMullen \cite{McM87} which implies that a nest of measurable sets in an Euclidean space 
shrinks to a set of full Hausdorff dimension. 
The criterion is employed in Section~\ref{sec:H-dim-PC} for the post-critical sets.
For the set of hairs without the end points, we directly employ the definition of the Hausdorff dimension.
%This criterion is also used in \cite{McM87} in order to study the Lebesgue measure and Hausdorff dimension of the Julia 
%sets of some transcendental entire functions. 
%Below we present the criterion.

Given measurable sets $K$ and $\Omega$ in $\mathbb{C}$, let $\area(K)$ denote the two-dimensional Lebesgue measure of $K$, and 
when $\area(\Omega)>0$, we define, and denote, the \textit{density} of $K$ in $\Omega$ by 
\[\dens(K,\Omega)=\frac{\area(K\cap \Omega)}{\area(\Omega)}.\]
%to denote the \textit{density} of $K$ in $\Omega$.

\begin{defi}
Let $\MK_n = \{K_{n,i}: 1\leq i\leq l_n\}$, for $n \geq 0$, be a finite collection of measurable subsets $K_{n,i} \subseteq \mathbb{C}$. 
We say that $\{\MK_n\}^{\infty}_{n=0}$ satisfies the \textbf{nesting conditions} if for all $n\geq 0$,  

\begin{enumerate}
\item $l_0=1$, $K_{0, 1}$ is a bounded and connected subset of $\mathbb{C}$, every $K_{n,i}$ has positive area,  
\item every $K_{n+1,i} \in \MK_{n+1}$ is contained in a (unique) $K_{n,j}\in \MK_n$, %where $1\leq i\leq l_{n+1}$ and $1\leq j\leq l_n$;
\item every $K_{n,i} \in \MK_n$ contains a $K_{n+1,j}\in \MK_{n+1}$, and %where $1\leq i\leq l_n$ and $1\leq j\leq l_{n+1}$;
\item $\area(K_{n,i}\cap K_{n,j})=0$, for all $1\leq i <j\leq l_n$. 
\end{enumerate}
 \end{defi}

Given a collection $\mathcal{F}$ of subsets of $\mathbb{C}$, we use the notation 
\[\hat{\mathcal{F}}:= \textstyle{\bigcup_{F\in \mathcal{F}}} F.\]

\begin{prop}
\label{prop:McMullen}
Assume that $\{\MK_n\}^{\infty}_{n=0}$ satisfies the nesting conditions, and there are sequences of positive numbers 
$(\delta_n)_{n\geq 0}$ and $(d_n)_{n\geq 0}$, with $d_n \to 0$ as $n\to \infty$, such that for $n\geq 0$ and $1\leq i \leq l_n$,
\[ \diam\, K_{n,i}\leq d_n, \qquad \text{ and } \qquad 
\dens \left(\hat{\MK}_{n+1}, K_{n,i}\right)\geq\delta_{n+1}.\]
Then,
\begin{equation}\label{equ:McM-form}
\dim_H\Big(\textstyle{\bigcap_{n \geq 0}} \hat{\MK}_n \Big)
\geq 2-\limsup_{n\to\infty}\frac{\sum_{k=1}^{n+1}|\log\delta_k|}{|\log d_n|}.
\end{equation}
\end{prop}

\begin{rmk}
Proposition~\ref{prop:McMullen} above and Proposition 2.2 in \cite{McM87} appear similar, 
but there is a minor difference.
%Our nest starts with $\MK_0$ while McMullen's begins with $\MK_1$. 
It seems that the upper bound in the sum in the statement of the latter proposition is $k$ rather than $k+1$. 
This difference is not crucial when dealing with the exponential maps, but is crucial in our case. 
For this reason, and for the reader's convenience, we present a proof below. 
\end{rmk}

\begin{proof}
By a rescaling, we may assume that $\area(\hat{\MK}_0)=1$. 
We inductively define the probability measure $\mu_n$ on $\hat{\MK}_n$, for $n\geq 0$.  
Let $\mu_0$ be the restriction of the two-dimensional Lebesgue measure on $\hat{\MK}_0$. 
Assuming $\mu_n$ is defined for some $n\geq 0$, $\mu_{n+1}$ is the unique measure such that on each $K_{n,j} \in \MK_n$, 
$\mu_{n+1}$ is a constant multiple of the Lebesgue measure satisfying $\mu_{n+1}(\hat{\MK}_{n+1} \cap K_{n,j})= \mu_n(K_{n,j})$. 

The sequence $\mu_n$ forms a martingale with the filtration $\{\MK_n\}_{n\geq 0}$.
Let $\mu$ be a (indeed, the) weak limit of $\mu_n$, as $n\to \infty$, 
which is a probability measure supported on $\hat{\MK} = \bigcap_{n \geq 0} \hat{\MK}_n$. 
It follows that for each $K_{n,j}$, we have $\mu(K_{n,j}) = \mu_n(K_{n,j}) \leq \mu_0(K_{n,j})/(\delta_n \delta_{n-1} \dots \delta_1)$. 

By Frostman's lemma, %(see for instance \cite[Theorem 8.8, p.\,112]{Mat95}), 
$\dim_H \hat{\MK} \geq s$ if there is a number 
$C(s)$ such that for all $a \in \C$ and $r>0$, $\mu(B_r(a))\leq C(s) r^s$. 
One only needs to consider this for small enough values of $r>0$.
Without loss of generality, we assume that $d_{n+1}<d_n$, for $n \geq 0$.

Choose $n\geq 0$ such that $d_{n+1} \leq r <d_n$, and let $\ML_{n+1}$ be the collection of all $K_{n+1,i}\in\MK_{n+1}$ which
meet $B_r(a)$. Then, $\hat{\ML}_{n+1} \subset B_{2r}(a)$, and we have
\[\mu(B_r(a)) \leq \mu(\hat{\ML}_{n+1}) \leq \frac{\mu_0(\hat{\ML}_{n+1})}{\delta_1\delta_2\cdots\delta_{n+1}}
\leq \frac{\pi 4 r^2}{\delta_1\delta_2\cdots\delta_{n+1}} 
\leq \pi 4 r^s\cdot\frac{d_n^{2-s}}{\delta_1\delta_2\cdots\delta_{n+1}}.\]
Define $b_n=d_n^{2-s}/(\delta_1\delta_2\cdots\delta_{n+1})$, for $n\geq 0$. 
If $s$ is smaller than the quantity on the right hand side of Equation~\eqref{equ:McM-form}, then we have 
$\limsup_{n\to\infty}b_n\leq 1$, and hence $(b_n)_{n \geq 0}$ is uniformly bounded from above. 
This means that $\hat{\MK}$ has Hausdorff dimension at least $s$.
\end{proof}

%\begin{rmk}
%If the diameter of each $K_{n,i}$ tends to zero much faster than the product of the densities
%$\delta_1\delta_2\cdots\delta_{n+1}$, as $n\to\infty$, then the superior limit in Equation~\eqref{equ:McM-form} 
%will be equal to zero and the Hausdorff dimension of $\bigcap_{n \geq 1}\MK_n$ will be equal to $2$.
%\end{rmk}
\section{Near-parabolic renormalisation scheme}\label{sec-renormalisation}
In Sections~\ref{subsec-IS-class} and \ref{subsec:near-para}, we present the Inou-Shishikura class of maps 
and the near-parabolic renormalisation operator. 
In Section \ref{subsec:renorm-tower} we define the renormalisation tower of a map, and in 
Section~\ref{subsec:chi-n-esti} we collect some estimates on the changes of coordinates 
which are used in this paper. 

\subsection{Inou-Shishikura class of maps}\label{subsec-IS-class}
Consider the cubic polynomial $P(z)=z(1+z)^2$. It has a parabolic fixed point at $0$ with multiplier $1$, and 
a critical point at $\cp_{P}=-1/3$ which is mapped to the critical value $\cv_P=-4/27$. 
There is another critical point at $-1$, which is mapped to $0$. Consider the ellipse 
\[E=\left \{x+ i y \in \C : ((x+0.18)/1.24)^2 + (y/1.04)^2 \leq 1\right\}, \]
and define 
\[U=\psi_1(\EC\setminus E), \text{~where~} \psi_1(z)=-4z/(1+z)^2.\]
The domain $U$ is a simply connected neighbourhood of $0$, which is symmetric with respect to the real axis, and contain $\cp_P$ but not $-1$. 
%, and $\overline{U}\cap(-\infty,-1]=\emptyset$

Following \cite[Section 4]{IS06}, we define 
\[\IS_0=
\left\{P \circ\varphi^{-1}: \varphi(U) \to \C
\mid 
\varphi:U \to \mathbb{C} \text{ is conformal}, ~\varphi(0)=0 \text{ and } \varphi'(0)=1 
\right\}.\]
Every map in $\IS_0$ has a parabolic fixed point at $0$, a unique critical point at $\varphi(-1/3)$ and a unique
critical value at 
\[\cv=-4/27. \]
For $\alpha \in \R$, we define
\[\IS_\alpha=\{f_0 (e^{2\pi i \alpha}z) \mid f_0 \in\IS_0\}.\]
For convenience, we also normalise the quadratic polynomials into the form 
\[Q_\alpha(z)= Q_0(e^{2\pi i \alpha}z)  = e^{2\pi i \alpha}z+ (27/16)e^{4\pi i \alpha}z^2,\]
so that all $Q_\alpha$ have a critical value at $-4/27$. 
%In particular, $Q_\alpha=Q_0\circ R_\alpha$, where $R_\alpha(z)=e^{2\pi i \alpha}z$.

Let us use the notation 
\[\QIS_\alpha = \{Q_\alpha\} \cup \IS_\alpha.\]

The class of maps $\IS_\alpha$ is in a one-to-one correspondence with the class of normalised univalent 
maps on the unit disk in $\mathbb{C}$, and hence by the Koebe distortion theorem, $\IS_\alpha$ is a
compact class of maps in the compact-open topology. 
Indeed, for all $f\in \QIS_0$, $f''(0)\neq 0$, and hence the parabolic fixed point of $f$ at $0$ is simple. 
It follows that if $\alpha \neq 0$ is small enough, the parabolic fixed point of $f_0$ at $0$ splits into 
two fixed points for the map $f(z)=f_0(e^{2\pi i \alpha}z)$, one at $0$ and another one at $\sigma_f$ near $0$. 
The fixed point $\sigma_f$ depends continuously on $f$.

\begin{prop}\label{P:fatou-coord}
There exist an integer $\kc \geq 1$ and a constant $r_1\in(0,1/2)$ satisfying $1/r_1-\textbf{k}\geq 2$ such that
for all $f \in \QIS_\alpha$ with $\alpha \in (0,r_1]$, there exist a univalent map $\Phi_f:\MP_f \to \C$ satisfying the following properties: 
\begin{enumerate}
\item $\MP_f$ is a Jordan domain bounded by piece-wise analytic curves, $\overline{\MP_f} \subset \Dom(f)$, 
and $0, \sigma_f, \cp_f \in \partial\MP_f$;
\item $\Phi_f(\cv)=1$ and 
\[\Phi_f(\MP_f)=\{\zeta\in\C:0<\re \zeta<1/\alpha-\textbf{k}\}\]
with $\im\Phi_f(z) \to +\infty$ as $z\to 0$ in $\MP_f$ and $\im\Phi_f(z)\to -\infty$ as $z\to \sigma_f$ in $\MP_f$;
\item $\Phi_f$ satisfies the Abel equation $\Phi_f(f(z))=\Phi_f(z)+1$ whenever $z,f(z) \in \MP_f$; 
\item $\Phi_f$ is unique and depends continuously on $f$.
\end{enumerate}
\end{prop}

\begin{figure}[!htpb]
  \setlength{\unitlength}{1mm}
  \centering
  \includegraphics[width=0.8\textwidth]{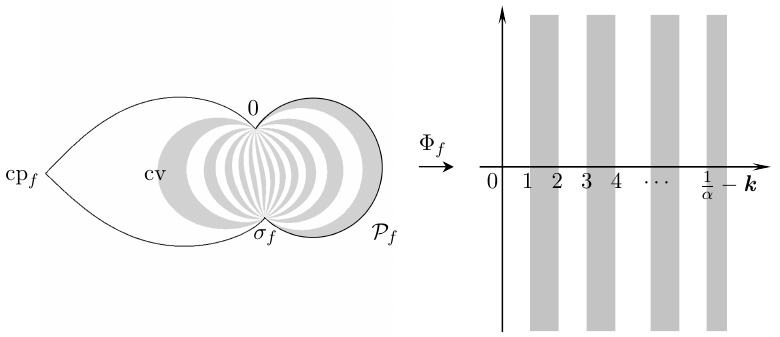}
  \caption{Illustration of Proposition~\ref{P:fatou-coord}.}
  \label{Fig:Petal-Phi}
\end{figure}

The map $\Phi_f$ in Proposition~\ref{P:fatou-coord} is called the (\emph{perturbed}) \textit{Fatou coordinate} of $f$, 
and $\MP_f$ is called a (\emph{perturbed}) \textit{petal} for $f$.
All parts of the above proposition, except for the existence of a uniform $\kc$, are proved in \cite[Main Theorem 1 and 3]{IS06}. 
The existence of a uniform $\kc$ was proved in \cite{Che09}, see also \cite[Proposition 12]{BC12} and \cite[Proposition 2.4]{Che19}. 
See Figure~\ref{Fig:Petal-Phi}.

\subsection{Near-parabolic renormalisation}\label{subsec:near-para}
\begin{figure}[!htpb]
  \setlength{\unitlength}{1mm}
  \centering
  \includegraphics[width=0.8\textwidth]{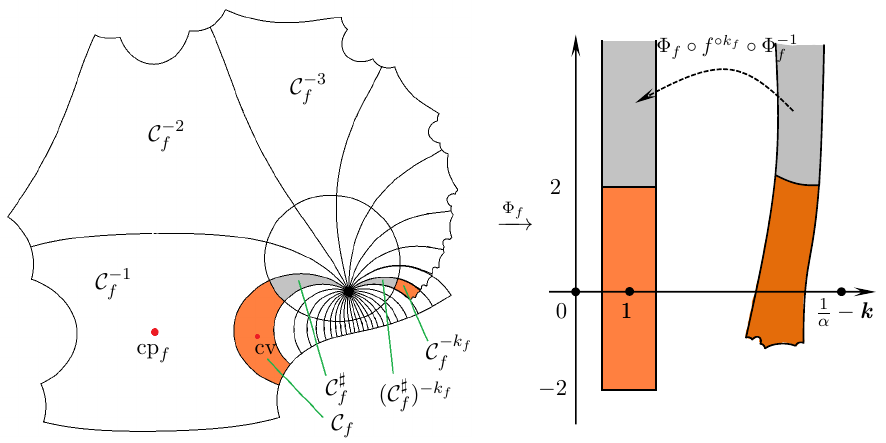}
  \caption{Illustration of the near-parabolic renormalisation.}
  \label{Fig:Near-parabolic}
\end{figure}

By Proposition~\ref{P:fatou-coord}, for $f \in \cup_{\alpha \in (0, r_1]} \QIS_\alpha$, we may define 
\[
\begin{split}
\MC_f=&\,\{z\in\MP_f:1/2\leq\re\Phi_f(z)\leq 3/2 \text{~and~} -2<\im\Phi_f(z)\leq 2\}, \\
\MC_f^\sharp=&\,\{z\in\MP_f:1/2\leq\re\Phi_f(z)\leq 3/2 \text{~and~} 2\leq\im\Phi_f(z)\}.
\end{split}
\]
See Figure~\ref{Fig:Near-parabolic}. 
Note that $\cv=-4/27\in \Int\, \MC_f$ and $0\in\partial \MC_f^\sharp$. 
Assume for a moment that there exists an integer $k_f \geq 1$, depending on $f$, satisfying the following properties:
\begin{enumerate}
\item For all $1\leq k\leq k_f$, there is a unique component $(\MC_f^\sharp)^{-k}$ of $f^{-k}(\MC_f^\sharp)$ containing
$0$ in its closure such that $f^{\circ k}:(\MC_f^\sharp)^{-k}\to\MC_f^\sharp$ is a univalent map;
\item There is a unique component $\MC_f^{-k}$ of $f^{-k}(\MC_f)$ intersecting $(\MC_f^\sharp)^{-k}$ such that
$f^{\circ k}:\MC_f^{-k}\to\MC_f$ is a covering of degree two ramified above $\cv$.
\item $\MC_f^{-k_f}\cup(\MC_f^\sharp)^{-k_f}$ is contained in $\{z\in\MP_f:1/2<\re\Phi_f(z)<\alpha^{-1}-\kc-1/2\}$.
\end{enumerate}
% Added needed property of compactly contained
Moreover, for all $ k = 1$, $\cdots$, $k_f $, the set $ ( \mathcal{C}_{ f } )^{ -k } \cup ( \mathcal{C}_f^{ \sharp } )^{ -k } $ is
compactly contained in $\Dom(f)$.

Let $k_f$ (assuming it exists) be the \emph{smallest} positive integer satisfying the above properties, and let
\[S_f=\MC_f^{-k_f}\cup(\MC_f^\sharp)^{-k_f}.\]
Consider the map 
\[\Phi_f\circ f^{\circ k_f}\circ\Phi_f^{-1}:\Phi_f(S_f)\to\C.\]
By Proposition~\ref{P:fatou-coord}-(c), the above map commutes with the translation by $+1$, and hence, 
it projects via the covering map
\[\Expo(\zeta)=-4e^{2\pi i \zeta}/27\]
to a map $\MMR f$ defined on a punctured neighbourhood of $0$. 
One may see that $\MMR f$ extends across $0$, with $\MMR f(0)=0$ and $(\MMR f)'(0)=e^{-2\pi i /\alpha}$. 
The map $\MMR f$ is called the \emph{near-parabolic renormalisation} of $f$.

\begin{thm}[{\cite[Main Theorem 3]{IS06}}]\label{thm-IS-attr-rep-3}
For all $\alpha \in (0, r_1]$ and $f \in \QIS_\alpha$, $\MMR (f) \in \{-1/\alpha\} \ltimes \MF$. 
%$\MMR f (z) = P \circ \psi^{-1}(e^{-2\pi i /\alpha }z)$, for some univalent map $\psi:U \to \C$ satisfying $\psi'(0)=0$ and $\psi'(0)=1$. 
%Moreover, $\psi$ extends to a univalent map from $U'$ to $\C$.
\end{thm}

Theorem~\ref{thm-IS-attr-rep-3} guarantees the existence of $k_f$ satisfying the properties in 
Section~\ref{subsec:near-para}. 

\begin{rmk}
A near-parabolic renormalisation scheme has been developed for a class of uni-critical maps in \cite{C14}, 
and independently for a class of uni-critical maps with local degree 3 in \cite{Yan15}. 
It is likely that the arguments in this paper may also be applied in those settings. 
\end{rmk}

\begin{prop}[{\cite[Prop.~2.7]{Che19}}]\label{prop-uniform-k-f}
There is $\kc_1 \in \mathbb{Z}$ such that for all $f \in \cup_{\alpha \in (0, r_1]} \QIS_\alpha$, $k_f\leq \kc_1$.
\end{prop}

%For another proof of Proposition~\ref{prop-uniform-k-f}, see \cite[Proposition 13]{BC12}. For the corresponding statements
%of Propositions \ref{P:fatou-coord} and \ref{prop-uniform-k-f} with $\alpha\in\C$ (specifically, $|\arg\alpha|<\pi/4$ and
%$|\alpha|$ is small), see \cite[Section 2]{ChSh14}.

\subsection{Renormalisation tower}\label{subsec:renorm-tower}
Let $s(z)= \overline{z}$ denote the complex conjugation map.
For any $f$ in $\cup_{\alpha \in [-r_1, 0)} \QIS_\alpha$, the map $\widetilde{f}=s\circ f\circ s$ satisfies 
$\widetilde{f}(0)=0$ and $(\widetilde{f})'(0)=e^{2\pi i (-\alpha)}$.
Because $s \circ f \circ s \in \QIS_0$, for all $f \in \QIS_0$ ($U$ and $P$ are real symmetric), we may extend the domain of definition of
$\MMR$ onto $\cup_{\alpha \in [-r_1, 0)} \QIS_\alpha$ by letting $\MMR(f)= \MMR(s \circ f \circ s)$. 

Let $r_1$ be the constant from Theorem~\ref{thm-IS-attr-rep-3}, and choose an integer $N$ satisfying 
\begin{equation}\label{E:N-condition}
N \geq 1/r_1+1/2,
\end{equation}
and let 
\[\HT_N =\left \{\varepsilon_0/(a_{0}+\varepsilon_1/(a_1+\varepsilon_2/(a_2+\dots ))) \mid \forall i\geq 0,  a_i \geq N, 
\varepsilon_i= \pm 1\right\}.\]
Because $1/\alpha_n = a_n + \varepsilon_{ n + 1 } \alpha_{ n + 1 }$, for all $n \geq 0$, one has 
\[\alpha_n^{-1}=a_{n}+\varepsilon_{n+1}\alpha_{n+1}\geq N-1/2\geq 1/r_1.\]

Fix $\alpha \in \HT_N$, and $f \in \QIS_\alpha$. We may define
\[f_0=
\begin{cases}
f  & \text{ if } \varepsilon_0=+1,\\
s\circ f\circ s & \text{ if } \varepsilon_0=-1.
\end{cases}\]
Successive applications of Theorem~\ref{thm-IS-attr-rep-3}, for $n\geq 1$, produces the following sequence of maps 
\[f_n=
\begin{cases}
\MMR(f_{n-1}) & \text{ if } \varepsilon_n=-1,\\
s\circ \MMR(f_{n-1})\circ s & \text{ if } \varepsilon_n=+1.
\end{cases}\]
Let $U_n=\Dom(f_n)$ for $n\geq 0$. Then for all $n\geq 1$, we have
\[f_n\in\IS_{\alpha_n},~f_n:U_n\to \C, ~f_n(0)=0,~f_n'(0)=e^{2\pi i \alpha_n} \text{~and~} \cv_{f_n}=-4/27.\]
For $n\geq 0$, let $\Phi_n = \Phi_{f_n}: \MP_n \to \mathbb{C}$ denote the Fatou coordinate of $f_n$, where 
$\MP_n=\MP_{f_n}$ is the Fatou petal of $f_n$. 
Then, for $n\geq 0$, $\MC_n=\MC_{f_n}$ and $\MC_n^\sharp=\MC_{f_n}^\sharp$ be the corresponding 
sets for $f_n$ defined in Section~\ref{subsec:near-para}.
Let $k_n=k_{f_n}$ be the smallest positive integer that appears in the definition of $\MMR$ such that
\[S_n=\MC_n^{-k_n}\cup(\MC_n^\sharp)^{-k_n}\subset\{z\in\MP_n:1/2<\re\Phi_n(z)<\alpha_n^{-1}-\kc-1/2\}.\]
Let $\sigma_n=\sigma_{f_n} \in \partial \MP_n$. 
By the compactness of $\IS$, $|\sigma_n|/\alpha_n$ is uniformly bounded from above and below.

\subsection{Successive changes of coordinates}\label{subsec:change}
For $n\geq 0$, define 
\[\Pi_n=\{\zeta\in\C:-1/(2\alpha_n)\leq\re\zeta<0, \,\im\zeta>0\} 
\bigcup\Phi_n(\MP_n) \bigcup \textstyle{\bigcup_{i=0}^{k_n+\gauss{1/(2\alpha_n)}}}(\Phi_n(S_n)+i).\]
Using $\Phi_n^{-1}(\zeta+1)=f_n\circ\Phi_n^{-1}(\zeta)$, we may extend 
$\Phi_n^{-1}:\Phi_n(\MP_n)\to\MP_n$ to a holomorphic map
\[\Phi_n^{-1}:\Pi_n\to U_n\setminus\{0\}. \]
Then, $\Phi_n^{-1}:\Pi_n\to\C\setminus\{0\}$ and $s\circ \Phi_n^{-1}:\Pi_n\to\C\setminus\{0\}$ may be lifted via the covering map $\Expo:\C\to\C\setminus\{0\}$ to define the holomorphic or an anti-holomorphic maps 
$\chi_n:\Pi_n\to\C$ according to 
\[\begin{cases}
\Expo\circ\chi_n = \Phi_n^{-1} & \text{ if } \varepsilon_n=-1,\\
\Expo\circ\chi_n = s\circ\Phi_n^{-1}  & \text{ if } \varepsilon_n=+1.
\end{cases}\]
Because $\Expo(1)=-4/27$ and $\Phi_n(-4/27)=1$, we may choose the branch of the lift so that 
$\chi_n(1)=1$, for all $n\geq 0$. 
Then, for each $j \in \Z$, define 
\[\chi_{n,j}=\chi_n+j.\]

\subsection{Some estimates on the changes of coordinates}\label{subsec:chi-n-esti}
The basic idea for controlling the changes of coordinates $\chi_n$ is to compare $\Phi_n^{-1}$ to the covering 
map $\tau_n(w)=\sigma_n/(1-e^{-2\pi i \alpha_n w})$ from $\mathbb{C}$ to the Riemann sphere minus 
$\{0, \sigma_n\}$, see for instance \cite{Sh00}. 
That is, one decomposes $\Phi_n^{-1}$ as $\tau_n \circ  L_n$, and aims to show that $L_n$ is close to the 
identity map. 
This has been extensively studied in \cite{Che13} and \cite[section 6]{Che19}. 
Below, we summarise some of those estimates that will be employed in this paper. 
See also \cite{AC18,CC15,FSh18} for some variations of those estimates and their applications.

\begin{prop}\label{P:est-imag}
For all $D_0>0$, there are $M_0, \widetilde{M}_0>0$ such that for all $n\geq 1$, we have 
\begin{enumerate}
\item If $\zeta \in \Pi_n$ with $\im\zeta\geq D_0/\alpha_n$, then
\[\Big|\im \chi_n(\zeta)-\Big(\alpha_n\,\im\zeta+\frac{1}{2\pi}\log\frac{1}{\alpha_n}\Big)\Big|\leq M_0.\]
\item If $\zeta\in\Pi_n$ with $\im \zeta\in[-2,D_0/\alpha_n]$, then
\[\Big|\im \chi_n(\zeta)-\tfrac{1}{2\pi}\min\big\{\log(1+|\zeta|),\log(1+|\zeta-1/\alpha_n|)\big\}\Big|\leq \widetilde{M}_0.\]
\end{enumerate}
\end{prop}

\begin{prop}\label{prop:esti-chi}
For all $D_0>0$, there are $M_1, \widetilde{M}_1\geq 1$ such that for all $n\geq 1$, we have
\begin{enumerate}
\item If $\zeta\in\Pi_n$ with $\im \zeta\geq -2$, $|\zeta|\geq D_0/\alpha_n$ and
$|\zeta-1/\alpha_n|\geq D_0/\alpha_n$, then
\[|\chi_n'(\zeta)-\alpha_n|\leq M_1\alpha_n e^{-2\pi\alpha_n \im\zeta}.\]
\item If $\zeta \in \Pi_n$ with $\im\zeta\geq -2$, and $1\leq |\zeta|<D_0/\alpha_n$ or
$1\leq |\zeta-1/\alpha_n|<D_0/\alpha_n$, then 
\[\widetilde{M}_1^{-1}\leq \min\{|\zeta|,|\zeta-1/\alpha_n|\}\cdot|\chi_n'(\zeta)|\leq \widetilde{M}_1.\]
\end{enumerate}
\end{prop}

In Propositions \ref{P:est-imag} and \ref{prop:esti-chi} one may make $M_0$ and $M_1$ arbitrarily small 
by making $D_0$ large enough. 
\section{The nest of partitions}\label{sec-partition}
In this section we build the nest of partitions for $\Lambda(f)$ to which we will apply the criterion 
in Section~\ref{sec-criterion}. 
Our approach mainly follows the construction of puzzle pieces for $\Lambda(f)$ that are used for the proof of 
Theorem~\ref{thm-Lambda-topo}.  
We do not need the full collection of the pieces from that paper, since we do not need all of $\Lambda(f)$ 
in order to achieve dimension $2$. 
The partition pieces are built at deep levels of renormalisation and mapped to the original phase space via the maps $\chi_n$. 
%There are many ways to build such partition pieces, however it will be technically simpler if we employ a 
%choice of partitions introduced in \cite{Che17} which enjoy equivariant properties with respect to the changes of 
%coordinates $\chi_n$. 
It is more convenient to work in Fatou coordinates because of the simpler behaviour of $\chi_n$ there. 
Later in this section we employ the estimates of Section~\ref{subsec:chi-n-esti} on $\chi_n$ to control the relative 
densities and diameters of those partition pieces. 

\subsection{Critical value arc}\label{subsec:topo}
By definition, a \textit{Cantor bouquet} is a compact subset of the plane which is homeomorphic to a set of the form 
\[\{re^{2 \pi i \theta} \in \C \mid \forall \theta \in \mathbb{R}/\mathbb{Z}, 0 \leq r \leq R(\theta) \}\]
where $R: \R/\Z \to [0, \infty)$ satisfies
\begin{itemize}
\item[(a)] both sets $R^{-1}(0)$ and $(\R/\Z) \setminus R^{-1}(0)$ are dense in $\R/\Z$, 
\item[(b)] for every $\theta_0\in \R/\Z$, $\limsup_{\theta \to \theta_0^+ } R(\theta) = R(\theta_0) = \limsup_{\theta \to \theta_0^-} R(\theta)$.
\end{itemize}

A \textit{one-sided hairy Jordan curve} is a compact subset of the plane which is homeomorphic to a set of the form 
\[\{re^{2\pi i \theta} \in \C \mid \forall \theta \in \R/\Z, 1 \leq r \leq 1+ R(\theta) \}\]
where $R: \R/\Z \to [1, \infty)$ satisfies (a) and (b) above. 
%\begin{itemize}
%\item[(a)] $R^{-1}(1)$ is dense in $\R/\Z$,  
%\item[(b)] $(\R/\Z) \setminus R^{-1}(1)$ is dense in $\R/\Z$, 
%\item[(c)] for each $\theta_0\in \R/\Z$ we have 
%\[\limsup_{\theta \to \theta_0^+ } R(\theta) = R(\theta_0) = \limsup_{\theta \to \theta_0^-} R(\theta).\]
%\end{itemize}

\begin{thm}[\cite{Che17}]\label{thm-Lambda-topo}
For every $\alpha \in \HT_N$ and every $f \in\QIS_\alpha$, the following hold: 
\begin{itemize}
\item[(i)] if $\alpha \in \HH$, $\Lambda(f)$ is a Jordan curve, 
\item[(ii)] if $\alpha\in\HB \setminus\HH$, $\Lambda(f)$ is a one-sided hairy Jordan curve enclosing $0$, 
%and the connected component of $\Lambda(f)\setminus\overline{\Delta}_f$ containing the critical value of $f$ is a $C^1$ curve;
\item[(iii)] if $\alpha \not \in \HB$, $\Lambda(f)$ is a Cantor bouquet at $0$. 
%and the connected component of $\Lambda(f) \setminus\{0\}$ containing the critical value of $f$ is a $C^1$ curve.
\end{itemize}
\end{thm}

Based on a toy model for the renormalisation, it is conjectured in \cite{Che23} that the above topological 
trichotomy holds for all rational functions. 
However, as illustrated by a toy model in a recent PhD thesis \cite{Russell26}, the geometry of the post-critical set 
appears drastically different when several recurrent critical points are interacting with the fixed point.

From now on, we assume that $\alpha \in \HT_N$ and $f\in \QIS_\alpha$. 

Let $\Delta(f)$ denote the Siegel disk of $f$ when $\alpha \in \HB$, and let $\Delta(f)=\{0\}$ 
when $\alpha \not \in\HB$. 
When $\alpha \in \HB$, $\Delta(f)$ is the region enclosed by the Jordan curve in $\Lambda(f)$.  
When $\alpha \notin \HH$, every connected component of $\Lambda(f) \setminus  \overline{\Delta(f)}$ is a Jordan arc, 
and the critical point of $f$ lies on an end point of such an arc.  

For $\alpha \notin \HB$, let $\Gamma_f$ be the connected component of $\Lambda(f) \setminus \{0\}$ containing
the critical value $-4/27$. 
Similarly, when $\alpha \in \HH$, let $\Gamma_f$ be an arc in $\overline{\Delta(f)}\setminus\{0\}$ which connects 
$0$ to $-4/27$. When $\alpha \in \HB \setminus \HH$, let $\Gamma_f$ be an arc in $\Lambda(f) \cup \Delta(f)\setminus \{0\}$ 
which connects $0$ to $-4/27$. In that case, $\Gamma_f$ is the union of the connecting component of 
$\Lambda(f) \setminus \overline{\Delta(f)}$ containing $-4/27$ and a ray in $\overline{\Delta(f)}\setminus\{0\}$ connecting $0$ 
to that component. 
In all these cases, we may choose $\Gamma_f$ so that its iterates are pairwise disjoint. 
Note that $\Gamma_f$ does not contain $0$ but it contains $-4/27$. 

\subsection{Partitions in a deep level}\label{subsec:partitions-deep}
\begin{figure}[!htpb]
  \setlength{\unitlength}{1mm}
  \centering
  \includegraphics[width=0.8\textwidth]{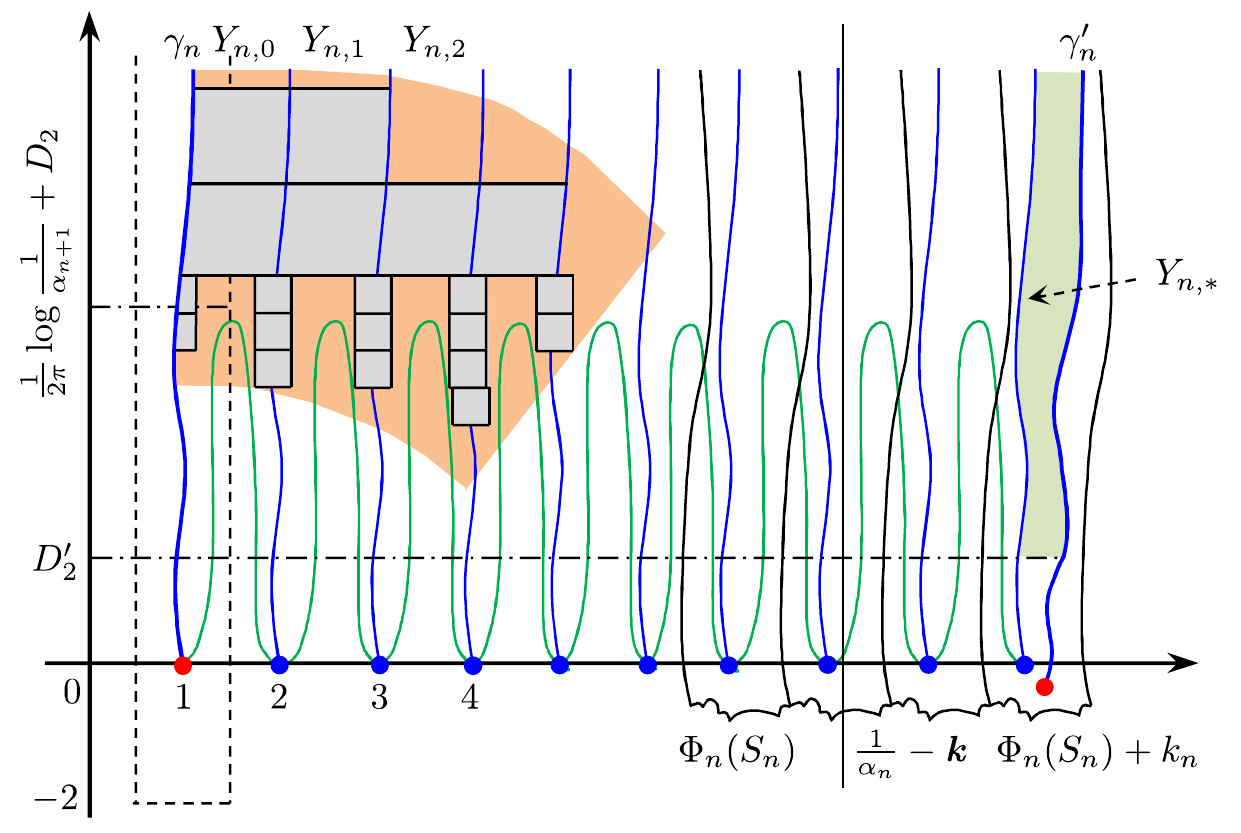}
  \caption{Illustration of the pieces on level $n$.}
  \label{Fig:Box}
\end{figure}

Here we extract some properties of $\Gamma_f$ which are used in this paper. 
See Sections 6.3, 6.4, 7.1, and 7.2 in \cite{Che17}, and Proposition 5.3 in \cite{FSh18}. 
See Figure~\ref{Fig:Box} for an illustration of the sets defined below.  

The set $\Gamma_f$ is contained in the Fatou petal $\MP_f$, and 
\[\gamma_f:=\Phi_f(\Gamma_f)\subset \mho:=\{\zeta \in \C \mid 1/2<\re\zeta<3/2, \im\zeta>-2\}.\]
By the definition of $\MMR(f)$ in Section~\ref{subsec:near-para}, there is a unique arc 
$\gamma_f'\subset \Phi_f(S_f)+k_f$ such that $\Phi_f^{-1}(\gamma_f')=\Gamma_f$.
The real parts of $\gamma_f$ and $\gamma_f'$ tend to finite limits as their imaginary parts
tend to $+\infty$, and they become straighter, in the sense that for $\zeta$ and $\zeta'$ in $\gamma_f$, 
or $\gamma_f'$, $\arg(\zeta-\zeta')$ can be made arbitrarily close to $\pm\pi/2$ 
by making $\im(\zeta)$ and $\im(\zeta')$ large enough, independent of $f$. 

For $n\geq 0$, let $\Gamma_n= \Gamma_{f_n}$, $\gamma_n= \gamma_{f_n}$ and $\gamma_n'= \gamma_{f_n}'$. 
For each $n\geq 0$, $\Phi_n^{-1}(\gamma_n)=\Gamma_n=\Phi_n^{-1}(\gamma_n')$, and $\chi_{n+1}(\gamma_{n+1})= \gamma_n$. 
Moreover, $\gamma_n+\mathbb{Z}$ and $\gamma_n'+\mathbb{Z}$ are disjoint sets, since 
they project to disjoint sets under $\Expo$. 

For $y \in \R$, let us define 
\[\BL_y=\{z \in \C \mid \im z=y\} \text{\quad and \quad} \mathbb{H}_y=\{z \in \C \mid \im z\geq y\}.\]
By the preceding paragraph, there is $D_2'>1$ such that for all $y \geq D_2'$, all $n\geq 0$ and all $j\in \Z$, both sets 
$\BL_y \cap (\gamma_n+j)$ and $\BL_y\cap (\gamma_n'+j)$ consist of single points. 

For $n \geq 0$, let 
\[Y_n=Y_n(D_2')\]
denote the closure of the region in $\BH_{D_2'}$ bounded by $\gamma_n \cup \gamma_n' \cup \BL_{D_2'}$, 
minus $\gamma_n'$. 
The set $Y_n$ is simply connected and uniformly close to a half-infinite strip of horizontal width $1/\alpha_n$. 
It contains $\gamma_n$ and part of $\BL_{D_2'}$, but not $\gamma_n'$.
In the same fashion, for $j\in \Z$, let
\[Y_{n,j}=Y_{n,j}(D_2')\]
denote the closure of the region in $\BH_{D_2'}$ bounded by $(\gamma_n+j) \cup (\gamma_n+j+1) \cup \BL_{D_2'}$, 
minus $\gamma_n+j+1$. 
Evidently, $Y_{n,j}=Y_{n,0}+j$. 

For $n \geq 0$, we define 
\[J_n = a_{n} + (\varepsilon_{ n + 1 } - 1)/ 2, \qquad \J_n = \Z \cap [0, J_n - 1],  \qquad \widetilde{\J}_n=\Z \cap [0, J_n].\]
and a half-infinite strip
\[Y_{n,*}=Y_n \setminus \textstyle{\bigcup_{j\in\J_n}} Y_{n,j}.\]
It follows from $1/\alpha_n = a_n + \varepsilon_{ n + 1 } \alpha_{ n + 1 }$ that $\overline{Y_{n,*}} \subset Y_{n,J_n}$.

Finally, for $n\geq 0$, let 
\[Y_{n,\diamond}=Y_{n,\diamond}(D_2')\]
denote the closure of the region in $\BH_{D_2'}$ bounded by $\gamma_n' \cup (\gamma_n'-1) \cup \BL_{D_2'}$,  
minus $\gamma_n'$. 
One may see that $Y_{n, \diamond}$ is contained in $\Phi_n(S_n)+\{k_n-1, k_n\}$.

%All of $Y_n$, $Y_{n,j}$, $Y_{n,*}$ and $Y_{n,\diamond}$ depend on the choice of the constant $D_2'\geq 1$. 

\subsection{Mapping up the tower}\label{subsec-go-tower}
Recall that $\overline{Y_n} \subset \Pi_n$; hence each $\chi_{n,j}$ is defined on $\overline{Y_n}$.

\begin{lemma}\label{lemma:tiled}
There is $D_2>0$ such that for all $n \geq 1$, $y_{n-1}=\tfrac{1}{2\pi}\log\tfrac{1}{\alpha_n}+D_2$ satisfies the following properties: 
\begin{enumerate}
\item If $\varepsilon_n=-1$, for all $j\in\J_{n-1}$, 
\[\chi_{n,j}( Y_n)\cap\BH_{y_{n-1}}= Y_{n-1,j}\cap\BH_{y_{n-1}}, \quad 
\chi_{n,J_{n-1}}( Y_n\setminus Y_{n,\diamond})\cap\BH_{y_{n-1}}= Y_{n-1,*}\cap\BH_{y_{n-1}}.\]
\item If $\varepsilon_n=+1$, for all $j\in\J_{n-1}$, 
\[\chi_{n,j+1}( \overline{Y_n} \setminus \gamma_n) \cap \BH_{y_{n-1}}= Y_{n-1,j}\cap\BH_{y_{n-1}}, \]
and
\[\chi_{n,J_{n-1}+1}( \overline{Y_{n,\diamond}} \setminus(\gamma_n'-1))\cap\BH_{y_{n-1}}
= Y_{n-1,*}\cap\BH_{y_{n-1}}.\]
\end{enumerate}
\end{lemma}

\begin{proof}
When $\varepsilon_n=-1$, $\chi_n:\Pi_n\to\C$ is holomorphic, $\chi_n(\gamma_n)=\gamma_{n-1}$, and 
$\chi_n(\gamma_n')=\gamma_{n-1}+1$.  
Then, the first relation in part (a) follows from Proposition~\ref{P:est-imag}. 

For the latter part of (a), by the definition of $\MMR$, 
$f_n(\Expo(\gamma_{n-1}'))=\Expo(\gamma_{n-1}+ J_{ n - 1 })=\Gamma_n$. 
This means that the set $\Gamma_n^{\cp}:= \Expo(\gamma_{n-1}')$ is an arc in $\Lambda(f_n) \cup \Delta(f_n)$
which connects the critical point of $f_n$ to $0$. That is, $\Gamma_n^{\cp}$ is the union of $\cp_n$ and the
component of $f_n^{-1}(\Gamma_n \setminus\{\cv\})$ with endpoints $0$ and $\cp_n$.
We also have $\Phi_n^{-1}(\gamma_n)=\Gamma_n$ and $\Phi_n^{-1}(\gamma_n'-1)=\Gamma_n^{\cp}$. 
%In particular, Proposition~\ref{P:est-imag} guarantees the existence of $D_2>0$ satisfying the latter relation in (a). 

When $\varepsilon_n=+1$, $\chi_n:\Pi_n\to\C$ is anti-holomorphic, $\chi_n(\gamma_n)=\gamma_{n-1}$, but 
$\chi_n(\gamma_n')=\gamma_{n-1}-1$. The rest of the proof is similar to the above. 
\end{proof}

\begin{rmk}
When $\varepsilon_n=-1$, the union of $\cup_{j \in \J_{n-1}} \chi_{n,j}(Y_n)$ and 
$\chi_{n,J_{n-1}}(Y_n\setminus Y_{n,\diamond})$ covers the upper end of $ Y_{n-1}$ because 
\[\left ( Y_{n-1,*} \cup \textstyle{\bigcup_{j \in \J_{n-1}}} Y_{n-1,j}\right)\cap\BH_{y_{n-1}}
= Y_{n-1}\cap\BH_{y_{n-1}}.\]
A similar statement holds when $\varepsilon_n=+1$.
\end{rmk}

In order to unify the presentation in the sequel, for $n\geq 1$ and $j \in \Z$, we use the notation 
\[\chi_{n,*}=\chi_{n,J_{n-1}} \text{\quad and\quad} \chi_{n,*+j}=\chi_{n,*}+j.\]
As before, ``$*$'' is a symbol, not equal to $J_{n-1}$ for $n\geq 1$. 
In particular, $Y_{n-1,*}$ and $Y_{n-1,J_{n-1}}$ are distinct, with $Y_{n-1,*}$ strictly contained in $Y_{n-1,J_{n-1}}$. 
Along the same lines, we define the symbolic addition 
\[\gamma_n+*=\gamma_n' \text{\quad and\quad} \gamma_n+(*-1)=(\gamma_n+*)-1=\gamma_n'-1, \quad \forall n\geq 0.\]
For $n\geq 1$, we define 
\[X_{n-1}=
 \left\{
\begin{aligned}
& \textstyle{\bigcup_{j\in\J_{n-1} \cup\{*\}}} \big(\chi_{n,j}( Y_n)\cap Y_{n-1,j}\big) & \text{ if } \varepsilon_n=-1,\\
& \textstyle{\bigcup_{j\in\J_{n-1}\cup\{*\}}} \big(\chi_{n,j+1}( \overline{Y_n})\cap Y_{n-1,j}) \big)  & \text{ if } 
\varepsilon_n=+1.
\end{aligned}
\right.\]
This is a connected set. 

The restriction of each $\chi_{n,j}$ to $\overline{Y_n}$ is injective. 
For $n\geq 1$, define $\xi_n: X_{n-1}\to Y_n$ as follows 
\begin{itemize}
\item if $\varepsilon_n=-1$, for $\zeta\in \chi_{n,j}( Y_n)\cap Y_{n-1,j}$ with $j\in\J_{n-1}\cup\{*\}$, define
\[\xi_n(\zeta)=\chi_{n,j}^{-1}(\zeta)= \chi_{n,0}(\zeta+j), \]
\item if $\varepsilon_n=+1$, for $j\in\J_{n-1}\cup\{*\}$, define
\[\xi_n(\zeta)=
\left\{
\begin{aligned}
& \chi_{n,0}^{-1}(\zeta) & ~~~\text{if} &\quad \zeta\in \chi_{n,1}(\overline{Y}_n)\cap\gamma_{n-1},\\
& \chi_{n,j+1}^{-1}(\zeta)  & ~~~\text{if} &\quad \zeta\in \chi_{n,j+1}(\overline{Y}_n)\setminus(\gamma_{n-1}+j).
\end{aligned}
\right.\]
\end{itemize}
%In both cases $\varepsilon_n=\pm 1$, we have $\xi_n\circ\chi_{n,j}(\zeta)=\zeta$ for all $\zeta\in Y_n$ and all $j\in\J_{n-1}\cup\{*\}$.
The map $\xi_n: X_{n-1}\to Y_n$ is a periodic function of period $+1$, and is the inverse of each $\chi_{n,j}$. 
Evidently, this map is not continuous along the arcs $(\gamma_{n-1}+j)\cap X_{n-1}$, for $1\leq j\leq J_{n-1}$. 
%For example, $\zeta\in (\gamma_{n-1}+1)\cap X_{n-1}$  is a boundary point of $Y_{n-1,0}$ and is also a boundary point of $Y_{n-1,1}$. If $\varepsilon_n=-1$, then by definition we have $\xi_n(\zeta)\in\gamma_n$. But there exists a sequence $(\zeta_k)_{k\in\N}\subset Y_{n-1,0}$ which converges to $\zeta$ such that $\xi_n(\zeta_k)$ converges to a point on $\gamma_n'$ as $k\to\infty$.
%To avoid confusion, when we look the image $\xi_n(\zeta)$ for
%$\zeta\in (\gamma_{n-1}+j)\cap X_{n-1}$, we will emphasize that $\zeta\in Y_{n-1,j-1}$ or $Y_{n-1,j}$.
%The ambiguity of the definition of $\xi_n$ on $(\gamma_{n-1}+j)\cap X_{n-1}$ does not effect the calculation of the
%Hausdorff dimension later since each curve $\gamma_{n-1}+j$ is $C^1$ continuous (see Theorem~\ref{thm-Lambda-topo}).
We will use $(\chi_n)_{n\geq 1}$ and $(\xi_n)_{n\geq 1}$ to map the partition pieces up and down the renormalisation tower. 

For $\zeta_0\in\C$ and $r>0$, let 
\[\Boxx(\zeta_0,r)=\{\zeta\in\C:|\re(\zeta-\zeta_0)|\leq r\text{ and }|\im(\zeta-\zeta_0)|\leq r\}.\]
Given positive numbers $D_3$, $D_3'\geq D_2'$, $\nu\in(0,1/2)$ and integers $n \geq 0$, let 
\[\begin{aligned}
\Xi_n=&~\Xi_n(D_3',D_3,\nu)\\
=&~Y_n\big(\tfrac{1}{2\pi}\log\tfrac{1}{\alpha_{n+1}}+D_3\big) \bigcup \left(Y_n(D_3')\cap 
B_\nu \left ( \gamma_n+\left (\widetilde{\J}_n \cup \{*,*-1\}\right ) \right ) \right) \\
\Xi_{n,j}&=\Xi_n\cap Y_{n,j}(D_3'), \qquad \forall j\in\J_n\cup\{*\}.
\end{aligned}\]
Also, for $n \geq 0$, real numbers $D_3''' > D_3'' \geq D_3$, and $j \in \J_n \cup \{*\}$, let 
\[W_{n,j}=W_{n,j}(D_3'',D_3''')=
\{\zeta\in Y_{n,j}: D_3''\leq\im\zeta-\tfrac{1}{2\pi}\log\tfrac{1}{\alpha_{n+1}}\leq D_3'''\}.\]

\begin{lemma}\label{lemma:fixed-width}
There exist constants $D_3$, $D_3'\geq D_2'$ and $\nu_0\in(0,1/20]$ such that for all $n\geq 1$, we have
\begin{enumerate}
\item $\Xi_{n-1}=\Xi_{n-1}(D_3',D_3,\nu_0)\subset X_{n-1}$;
\item if $\zeta_{n-1}\in\gamma_{n-1}+j'$ and $\zeta_{n-1}'\in\Xi_{n-1,j}\cap B_{\nu_0}(\gamma_{n-1}+j')$, for some 
$j \in\J_{n-1}\cup\{*\}$ and $j'\in\widetilde{\J}_{n-1}\cup\{*,*-1\}$, and $\im\zeta_{n-1}'\geq \im\zeta_{n-1}-\nu_0$, then 
\[\im\xi_n(\zeta_{n-1}')\geq \tfrac{3}{4}\,\im\xi_n(\zeta_{n-1});\]
\item if $\zeta_{n-1}\in\gamma_{n-1}+j'$ satisfies $\Boxx(\zeta_{n-1},\nu_0)\cap Y_{n-1,j}\neq\emptyset$, for some 
$j \in \J_{n-1}\cup\{*\}$ and $j' \in \widetilde{\J}_{n-1}\cup\{*,*-1\}$, then $\xi_n:\Boxx(\zeta_{n-1},\nu_0)\cap Y_{n-1,j}\to Y_n$ has 
a univalent or anti-univalent extension to $\widetilde{\xi}_{n,j}:\Boxx(\zeta_{n-1},20\nu_0)\to\Pi_n$; 
\item for any $D_3'''>D_3''\geq D_3$ and $j\in\J_{n-1}\cup\{*\}$, $\xi_n:W_{n-1,j}\cap Y_{n-1}\to Y_n$ has a univalent or anti-univalent extension 
$\widetilde{\xi}_{n,j}:B_{\nu_0}(W_{n-1,j})\to\Pi_n$. 
\end{enumerate}
\end{lemma}

Part (a) of the above lemma plays a key role in estimating relative densities in later sections. 
Part (b) will be used to locate the positions of the boxes when mapping to deeper or shallower levels of the renormalisation tower.
When $\xi_n$ is defined on a ``half" box (e.g.\ when $\zeta_{n-1} \in \gamma_{n-1} +\{0, *\}$, we consider the left or the right ``half" 
part of the box, Parts (c) and (d) of Lemma~\ref{lemma:fixed-width} are used to control the distortion of $\xi_n$. 

\begin{proof}
We present the proof for $\varepsilon_n=-1$, which produces univalent extensions. 
The other case may be proved along the same lines, yielding anti-univalent extensions. 

(a) By Proposition~\ref{P:est-imag}, $Y_n((1/(2\pi) \log(1/\alpha_{n+1})+D_3)$ is contained in $X_n$, provided $D_3$ is large enough. 
Thus, we need to show that $Y_n(D_3')\cap B_\nu( \gamma_n+(\widetilde{\J}_n \cup \{*,*-1\} ))$
is contained in $X_n$. 
%Indeed, by the definition of $X_n$, it is enough to show that $Y_n(D_3')\cap B_\nu( \gamma_n+\{0, *\} )) \subset X_n$. 
 
Recall $D_2'\geq 1$ from Section~\ref{subsec:partitions-deep}. 
Because $\gamma_n \subset \mho$, by the pre-compactness of $\bigcup_{\alpha \in (0,r_1]}\QIS_\alpha$ and 
Proposition~\ref{prop-uniform-k-f}, there is $C_0'>0$ such that for all $n \geq 0$, (see \cite[Lemma~4.13]{Che17} for more details)  
\begin{equation}\label{equ:gamma-pri-posi}
|\re\zeta-1/\alpha_n|\leq C_0', \qquad \forall \zeta \in \gamma_n'.
\end{equation}
By Lemma~\ref{lemma:tiled}(a), $\chi_{n,J_{n-1}}(\gamma_n'-1)=\gamma_{n-1}'$. 
Both $\chi_n=\chi_{n,0}:\gamma_n\to\gamma_{n-1}$ and $\chi_n:\gamma_n'\to\gamma_{n-1}+1$ are homeomorphisms. 
Then, by Proposition~\ref{P:est-imag}, there are $C_1 \geq D_2'$ and $C_1'>0$ such that for all $n \geq 1$, 
\begin{itemize}
\item any $\zeta_{n-1}\in(\gamma_{n-1}+1)\cap\BH_{C_1}$ has a unique preimage $\zeta_n\in\gamma_n\cap\BH_{2D_2'}$ 
under $\chi_{n,1}$;
\item any $\zeta_{n-1}' \in(\gamma_{n-1}+1)\cap\BH_{C_1}$ has a unique preimage $\zeta_n'\in\gamma_n'\cap\BH_{2D_2'}$ 
under $\chi_n=\chi_{n,0}$;
and 
\item any $\zeta_{n-1}'' \in \gamma_{n-1}'\cap\BH_{C_1}$ has a unique preimage $\zeta_n''\in(\gamma_n'-1)\cap\BH_{2D_2'}$ 
under $\chi_{n,J_{n-1}}$;
%\item for all $\zeta\in\BL_{D_2'}\cap Y_n$, $\im\chi_n(\zeta)\leq \tfrac{1}{2\pi}\log\frac{1}{\alpha_n}+C_1'$.
\end{itemize}

We consider two cases depending on the location of $\zeta_n\in\gamma_n\cap\BH_{2D_2'}$. 
When $\im\zeta_n\leq 1/\alpha_n$, let  
\[V_n^+=\{\zeta\in Y_n: 3/4 \leq \im \zeta /\im \zeta_n \leq 4/3 \text{ and }\re\zeta\leq \im(\zeta_n/2)+1\}.\]
Recall that $\chi_{n,1}: Y_n\to Y_{n-1,1}$ has a univalent extension onto a neighbourhood of $\overline{Y}_n$. 
By Proposition~\ref{prop:esti-chi}(b), there is $\widetilde{M}_1\geq 1$ such that 
$\widetilde{M}_1^{-1}\leq |\chi_{n,1}'(\zeta)|/\im\zeta_n\leq \widetilde{M}_1$ for all $\zeta\in V_n^+$. 
This means that $\chi_{n,1}(V_n^+)$ is a topological disk satisfying 
\[|\zeta-\zeta_{n-1}|\geq \varrho_1 \text{\quad for all }\zeta\in \chi_{n,1}(\partial V_n^+\setminus\gamma_n),\]
where $0<\varrho_1<1$ is a constant depending only on $\widetilde{M}_1$.

The other case is when $\im \zeta_n> 1/\alpha_n$, where we let 
\[V_n^+=\{\zeta\in Y_n: -3/(4\alpha_n) \leq \im\zeta-\im\zeta_n \leq 3/(4\alpha_n), \re \zeta \leq 1/(2\alpha_n)+1\}.\]
In the same fashion, Proposition~\ref{prop:esti-chi}(a) implies that there is $\widetilde{\varrho}_1 \in (0,1)$, independent of $n$, 
such that the topological disk $\chi_{n,1}(V_n^+)$ satisfies 
\[|\zeta-\zeta_{n-1}| \geq \widetilde{\varrho}_1 \text{\quad for all }\zeta\in \chi_{n,1}(\partial V_n^+\setminus\gamma_n).\]

Similarly, $\chi_n: Y_n\to Y_{n-1,0}$ has a univalent extension onto a neighbourhood of $\overline{Y}_n$. 
When, $\im \zeta_n' \leq 1/\alpha_n$, let 
\[V_n^-=\{\zeta\in Y_n:3/4 \leq \im\zeta/\im\zeta_n' \leq 4/3, \re\zeta\geq 1/\alpha_n-\im\zeta_n/2-1\},\]
and when $\im\zeta_n> 1/\alpha_n$, consider
\[V_n^-=\{\zeta\in Y_n:-3/(4\alpha_n) \leq \im\zeta-\im\zeta_n' \leq 3/(4\alpha_n), \re\zeta\geq 1/(2\alpha_n)-1\}.\]
As in the above paragraphs, Proposition~\ref{prop:esti-chi} implies that there is $\varrho_2 \in (0,1)$, independent of $n$, 
such that 
\[|\zeta-\zeta_{n-1}|\geq \varrho_2 \text{\quad for all }\zeta\in \chi_n(\partial V_n^-\setminus\gamma_n').\]
Note that the interior of $\widetilde{V}_{n-1}=\chi_{n,1}(V_n^+)\cup\chi_n(V_n^-)$ is a neighbourhood of $\zeta_{n-1}$, and
for $\zeta\in\partial\widetilde{V}_{n-1}$,
\[|\zeta-\zeta_{n-1}|\geq \varrho'=\min\{\varrho_1,\widetilde{\varrho}_1,\varrho_2\}.\]
Then, if $D_3' \geq C_1+1$,  $Y_{n-1}(D_3') \cap B_{\varrho'}(\gamma_{n-1}+\Z) \subset X_{n-1}$.

By a similar argument, using \eqref{equ:gamma-pri-posi} and Proposition~\ref{prop:esti-chi}, there is $\varrho''>0$, 
independent of $n$, such that $Y_{n-1}(D_3')\cap B_{\varrho''}(\gamma_{n-1}')$ is contained in $X_{n-1}$.
%\begin{equation}
%Y_{n-1}(D_3')\cap B_{\varrho''}(\gamma_{n-1}')\subset X_{n-1}.
% \text{ and } B_{\varrho''}(\gamma_{n-1}')\cap\BH_{D_3'}\subset \textstyle{\bigcup_{j\in\N}} \chi_{n,j}(Y_n).
%\end{equation}
We set $\nu_0=\min\{\varrho', \varrho''\}/20$.  

\medskip

(b) Given $\zeta_{n-1}$ and $\zeta_{n-1}'$ as in part (b), by the definitions of $V_n^\pm$ and $\nu_0$ in Part (a), 
there exists $\widetilde{\zeta}_{n-1} \in \gamma_{n-1}+j'$ such that $\im\widetilde{\zeta}_{n-1}=\im\zeta_{n-1}'+\nu_0$ and hence 
\[\im\xi_n(\zeta_{n-1}')\geq 3\,\im\xi_n(\widetilde{\zeta}_{n-1})/4\geq 3\,\im\xi_n(\zeta_{n-1})/4.\]

\medskip

(c) %Let $\zeta\in(\gamma_{n-1}+\widetilde{\J}_{n-1}\cup\{*,*-1\})\cap X_{n-1}$. 
Suppose $\xi_n=\chi_{n,j}^{-1}$ is defined on $\Boxx(\zeta,\nu_0)\cap Y_{n-1,j}$, for some $j \in \J_{n-1}\cup\{*\}$. 
Because $\chi_n$ is defined on $\Pi_n$, and $V_n^\pm\subset \Pi_n$, the extensions are guaranteed by the choice of $\nu_0$. 

(d) First note that if $D_3'' > D_3$, then $W_{n-1,j}$ is contained in $X_{n-1}$. By adding $+1$ to $D_3''$, we may guarantee 
the existence of the univalent extension. 
\end{proof}

%We will use the following estimations, which can be seen as an inverse version of Proposition~\ref{P:est-imag} in some sense.

\begin{lemma}\label{lemma:chi-inverse}
For any $\varepsilon\in(0,1/10)$, there exist positive constants $D_4 \geq D_3$, $D_4' \geq D_3'$ and $\widetilde{M}_4 \geq 1$ 
such that for all $n \geq 1$, $\zeta_{n-1} \in \Xi_{n-1}(D_4',D_4,\nu_0)$ and $\zeta_n=\xi_n(\zeta_{n-1})$, 
we have
\begin{enumerate}
\item If $\im\zeta_{n-1} \geq \frac{1}{2\pi}\log\tfrac{1}{\alpha_n}+D_4$, then
\[\im\zeta_n\geq (16/9) \im\zeta_{n-1} \text{\quad and\quad} |\chi_n'(\zeta_n)-\alpha_n|\leq \alpha_n \varepsilon/10;\]
\item If $D_4'\leq \im\zeta_{n-1}\leq \frac{1}{2\pi}\log\tfrac{1}{\alpha_n}+D_4+2$, then
\[\im\zeta_n\geq (4/3) \im \zeta_{n-1} \text{\quad  and \quad}
\widetilde{M}_4^{-1} e^{-2\pi\im\zeta_{n-1}} \leq |\chi_n'(\zeta_n)|
\leq \widetilde{M}_4 e^{-2\pi\im\zeta_{n-1}} < 3/5. \]
\end{enumerate}
\end{lemma}

In Lemma~\ref{lemma:chi-inverse}, first we find $D_4\geq D_3$ so that Part (a) holds. We then choose $D_4'\geq D_3'$ so that 
Part (b) holds. If $D_4' > \frac{1}{2\pi}\log\tfrac{1}{\alpha_n}+D_4+2$ for some $n \geq 1$, the statement of Part (b) 
is void. 

\begin{proof}
(a) By Proposition~\ref{P:est-imag}(a), if $\im\zeta_n\geq D_0/\alpha_n> D_2'$ for some $D_0>0$, there is $M_0>0$ such that
\begin{equation}\label{equ:est-1}
\Big|\im \zeta_n-\frac{1}{\alpha_n}\Big(\im\zeta_{n-1}-\frac{1}{2\pi}\log\frac{1}{\alpha_n}\Big)\Big|\leq\frac{M_0}{\alpha_n}.
\end{equation}
If $\im\zeta_n<D_0/\alpha_n$, by Proposition~\ref{P:est-imag}(b), there is $M_0'>0$ such that 
$\im\zeta_{n-1}< \frac{1}{2\pi}\log\tfrac{1}{\alpha_n}+M_0'$.
Therefore, if $\im\zeta_{n-1}\geq\frac{1}{2\pi}\log\tfrac{1}{\alpha_n}+M_0'$, then $\im\zeta_n\geq D_0/\alpha_n> D_2'$
and \eqref{equ:est-1} holds.

Suppose that $\im\zeta_{n-1}\geq\frac{1}{2\pi}\log\tfrac{1}{\alpha_n}+M_0+M_0'$.
We denote $\im\zeta_{n-1}=\frac{1}{2\pi}\log\tfrac{1}{\alpha_n}+y$ with $y\geq M_0+M_0'$. 
Then by \eqref{equ:est-1} we have $\im\zeta_n\geq (y-M_0)/\alpha_n$ and
\[\frac{\im\zeta_n}{\im\zeta_{n-1}}\geq \frac{y-M_0}{\alpha_n y+\tfrac{1}{2\pi}\alpha_n\log\tfrac{1}{\alpha_n}}.\]
Note that $\tfrac{1}{2\pi}\alpha_n\log\tfrac{1}{\alpha_n}>0$ is uniformly bounded from above. Since $0<\alpha_n<1/2$, there
is a constant $M_0''>0$ such that for all $y \geq M_0''$, then $\im\zeta_n/\im\zeta_{n-1}\geq 16/9$.

On the other hand, if $\im\zeta_{n-1}\geq\frac{1}{2\pi}\log\tfrac{1}{\alpha_n}+M_0+M_0'$, we have $\im\zeta_n\geq M_0'/\alpha_n$.
By Proposition~\ref{prop:esti-chi}(a), there exists a constant $M_1\geq 1$ such that
\[|\chi_n'(\zeta_n)-\alpha_n|\leq M_1\alpha_n e^{-2\pi\alpha_n \im\zeta_n}.\]
If in addition 
\[\im\zeta_{n-1}\geq\frac{1}{2\pi}\log\tfrac{1}{\alpha_n}+M_0+M_0'+\tfrac{1}{2\pi}\log(10 M_1/\varepsilon),\] 
by \eqref{equ:est-1}, $\im \zeta_n \geq (1/(2\pi \alpha_n)) \log (10M_1/\varepsilon)$, and hence 
$|\chi_n'(\zeta_n)-\alpha_n|\leq \alpha_n\varepsilon/10$.
Therefore, Part (a) holds if we set $D_4=\max\{M_0'',M_0+M_0'+\tfrac{1}{2\pi}\log(10 M_1/\varepsilon)\}$.

\medskip

(b) Without loss of generality, we assume that $\varepsilon_n=-1$ and
$\zeta_{n-1}\in Y_{n-1,0}\cap\Xi_{n-1}(D_3',D_4,\nu_0)$. The rest is divided into two cases:
$\zeta_{n-1}\in B_{\nu_0}(\gamma_{n-1})\cap Y_{n-1,0}$ and $\zeta_{n-1}\in B_{\nu_0}(\gamma_{n-1}+1)\cap Y_{n-1,0}$.

Suppose that $\zeta_{n-1}\in B_{\nu_0}(\gamma_{n-1})\cap Y_{n-1,0}$. 
There exists $\zeta_{n-1}'\in\gamma_{n-1}$ with $\im\zeta_{n-1}'=\im\zeta_{n-1}$ such that 
$\zeta_n'=\xi_n(\zeta_{n-1}')=\chi_n^{-1}(\zeta_{n-1}')\in\gamma_n$. 
We also have $\zeta_n=\xi_n(\zeta_{n-1})=\chi_n^{-1}(\zeta_{n-1})$. 
Since $\im\zeta_{n-1}'\leq \frac{1}{2\pi}\log\tfrac{1}{\alpha_n}+D_4+2$, by Proposition~\ref{P:est-imag}(b), there exists a
constant $\widetilde{M}_0 >0$ depending on $D_4$ such that 
\[\big|\im \zeta_{n-1}'-\tfrac{1}{2\pi}\log(1+|\zeta_n'|)\big|\leq \widetilde{M}_0.\]
If $\im\zeta_{n-1}'\geq\widetilde{M}_0+1$, then we have
\[2\pi(\im\zeta_{n-1}'-\widetilde{M}_0)\leq\log(1+|\zeta_n'|)\leq 2\pi(\im\zeta_{n-1}'+\widetilde{M}_0).\]
Recall that $\zeta_n' \in \gamma_n$ is contained in $\mho=\{\zeta\in\C:1/2<\re\zeta<3/2 \text{ and }\im\zeta>-2\}$.
Then we have
\begin{equation}\label{equ:est-2}
C_1^{-1} e^{2\pi\im\zeta_{n-1}'}\leq\im\zeta_n'\leq C_1 e^{2\pi\im\zeta_{n-1}'},
\end{equation}
where $C_1=2e^{2\pi\widetilde{M}_0}$. 
Therefore, there is a constant $C_1'>0$ such that if $\im\zeta_{n-1}'\geq C_1'$,
then $\im\zeta_n'\geq \tfrac{16}{9}\,\im\zeta_{n-1}'$. By the definition of $\nu_0$ and Lemma~\ref{lemma:fixed-width}(b),
we have
\begin{equation}\label{equ:est-3}
\frac{3}{4}\,\im\zeta_n'\leq\im\zeta_n\leq \frac{4}{3}\,\im\zeta_n' \text{\quad and\quad}
\im\zeta_n\geq \frac{4}{3}\,\im\zeta_{n-1}.
\end{equation}

By Proposition~\ref{prop:esti-chi}(b), there is $\widetilde{M}_1\geq 1$ depending on $D_4$ 
such that $\widetilde{M}_1^{-1}/|\zeta_n|\leq|\chi_n'(\zeta_n)|\leq \widetilde{M}_1/|\zeta_n|$.
By \eqref{equ:est-2} and \eqref{equ:est-3}, that implies that there is $\widetilde{M}_4 \geq 1$ such that
\[\widetilde{M}_4^{-1} e^{- 2\pi\im\zeta_{n-1}} \leq |\chi_n'(\zeta_n)| \leq  \widetilde{M}_4 e^{-2\pi\im\zeta_{n-1}}.\]
Moreover, we assume that $C_1'>0$ is large enough such that if $\im\zeta_{n-1}\geq C_1'$,
then $\widetilde{M}_4/e^{2\pi\im\zeta_{n-1}}<3/5$.
Therefore, we may set $D_4'=\max\{D_3',\widetilde{M}_0+1,C_1'\}$. 
%then Part (b) holds under the assumption that $\zeta_{n-1}\in B_{\nu_0}(\gamma_{n-1})\cap Y_{n-1,0}$.

The argument is similar in the latter case $\zeta_{n-1}\in B_{\nu_0}(\gamma_{n-1}+1)\cap Y_{n-1,0}$, 
where one uses \eqref{equ:gamma-pri-posi}.
\end{proof}

For $\varepsilon\in(0,1/10)$, let $D_4=D_4(\varepsilon)$ and $D_4'=D_4'(\varepsilon)$ denote the corresponding constants in 
Lemma~\ref{lemma:chi-inverse}.  
For $n \geq 0$, define 
\begin{equation}\label{equ:h-n}
h_n=(4/3)^n D_4', 
\end{equation}
and 
\[T_n=T_n(\varepsilon)=\BH_{h_n}\cap\Xi_n(D_4',D_4,\nu_0).\]
In particular, $T_0=\Xi_0(D_4',D_4,\nu_0)$. 
By Lemma~\ref{lemma:fixed-width}-(a), $T_n\subset \Xi_n\subset X_n$ for all $n \geq 0$. 
Also, by Lemma~\ref{lemma:chi-inverse}, for each $n\geq 1$, 
\begin{equation}\label{equ:T-n-above}
\xi_n(T_{n-1})\subset \BH_{h_n}.
\end{equation}

Note that $D_4$ and $D_4'$ are positive numbers depending on $\varepsilon$, while $\nu_0\in(0,1/20]$ is a universal constant (independent on $\varepsilon$).
The following lemma will be used to estimate the diameters of some compact sets as we go up the renormalisation tower.

\begin{lemma}\label{lemma:diameter-n}
For any $\varepsilon\in(0,1/10)$ there is $M_2=M_2(\varepsilon)>0$ satisfying the following property.  
Assume that for some $ \zeta_0 \in T_0 =T_0(\varepsilon)$, we have $\zeta_n = \xi_n(\zeta_{n-1}) \in T_n$ for all $n \geq 1$. 
Then, for all $n \geq 1$, we have 
\[|\chi_n'(\zeta_n)|\leq \widetilde{\mu}_n<3/5,\]
where 
\[\widetilde{\mu}_n=\left\{
\begin{aligned}
& 11\alpha_n/10 & \text{ if } h_{n-1}\geq \tfrac{1}{2\pi}\log\tfrac{1}{\alpha_n}+D_4,\\
& M_2 e^{-2\pi h_{n-1}} & \text{ if } h_{n-1}< \tfrac{1}{2\pi}\log\tfrac{1}{\alpha_n}+D_4.
\end{aligned}
\right.\]
\end{lemma}

\begin{proof}
The case that $h_{n-1}\geq \tfrac{1}{2\pi}\log\tfrac{1}{\alpha_n}+D_4$ is immediate by Lemma~\ref{lemma:chi-inverse}(a).
If $h_{n-1}\leq \im\zeta_{n-1}<\tfrac{1}{2\pi}\log\tfrac{1}{\alpha_n}+D_4$, then by Lemma~\ref{lemma:chi-inverse}(b), we
have
\[|\chi_n'(\zeta_n)|\leq \widetilde{M}_4 e^{-2 \pi \im \zeta_{n-1}} \leq \widetilde{M}_4 e^{-2\pi h_{n-1}} < 3/5.\]
If $\im\zeta_{n-1}\geq\tfrac{1}{2\pi}\log\tfrac{1}{\alpha_n}+D_4>h_{n-1}=(\tfrac{4}{3})^{n-1} D_4'$, then by
Lemma~\ref{lemma:chi-inverse}(a), 
\[|\chi_n'(\zeta_n)|\leq \frac{11}{10} \alpha_n 
= \tfrac{11}{10}e^{2\pi D_4} e^{-2\pi(\tfrac{1}{2\pi}\log\tfrac{1}{\alpha_n}+D_4)}
\leq \tfrac{11}{10} e^{2\pi D_4} e^{- 2\pi h_{n-1}}.\]
Without loss of generality, we assume that $D_4' \geq D_4+1$, so that 
$(11/10) e^{2\pi D_4} e^{-2\pi h_{n-1}} < 3/5$. 
%We may set $M_2=\max\{\widetilde{M}_4,\tfrac{11}{10}e^{2\pi D_4}\}$. 
\end{proof}

\subsection{Boxes and almost rectangles}\label{subsec:box-rectangle}
Let us fix an arbitrary $\varepsilon\in(0,1/10)$, and let $D_4=D_4(\varepsilon)$ and $D_4'=D_4'(\varepsilon)$ 
be the corresponding constants from Lemma \ref{lemma:chi-inverse}. 
Without loss of generality, we may assume that $D_4'>D_2'+1$ is sufficiently large so that 
\begin{equation}\label{equ:small-arg}
|\arg(\zeta-\zeta')-\pi/2|<\arctan(\nu_0/5)\leq\arctan(1/100),
\end{equation}
where either $\zeta, \zeta' \in\gamma_n$ or $\zeta, \zeta' \in\gamma_n'$, and $\im\zeta>\im\zeta'\geq D_4'-1$.
For $n \geq 0$, we partition $Y_{n,J_n-1}$ as 
\[Y_{n,J_n-1}^- =Y_{n,*}-1 \quad \text{ and } \quad Y_{n,J_n-1}^+ = Y_{n,J_n-1}\setminus Y_{n,J_n-1}^-.\]

For $n \geq 0$ and $y=\tfrac{1}{2\pi}\log\tfrac{1}{\alpha_{n+1}}+D_4$, we say that $R \subset Y_n \cap \BH_y$ is 
an \textbf{almost rectangle} if there are real numbers $b > a \geq D_4$, with $1 \leq b-a \leq 3$, such that either 
\begin{itemize}
\item $R=W_{n,j}(a,b)$ for some $j \in \J_n \cup \{*\}$, or % and $b>a\geq D_4$ satisfying $1\leq b-a\leq 3$,
\item $R=W_{n, J_n-1} \cap  Y_{n,J_n-1}^-$, or %  \{\zeta \in Y_{n,J_n-1}^-: a \leq \im\zeta- \tfrac{1}{2\pi} \log\tfrac{1}{\alpha_{n+1}} \leq b\}$, 
%for some $b>a\geq D_4$ satisfying $1\leq b-a\leq 3$.
\item $R= W_{n, J_n-1} \cap  Y_{n,J_n-1}^+$.  % \{\zeta \in Y_{n,J_n-1}^+: a \leq \im\zeta- \tfrac{1}{2\pi} \log\tfrac{1}{\alpha_{n+1}} \leq b\}$. 
%for some $b>a\geq D_4$ satisfying $1\leq b-a\leq 3$.
\end{itemize}

\begin{figure}[!htpb]
  \setlength{\unitlength}{1mm}
  \centering
  \includegraphics[width=0.7\textwidth]{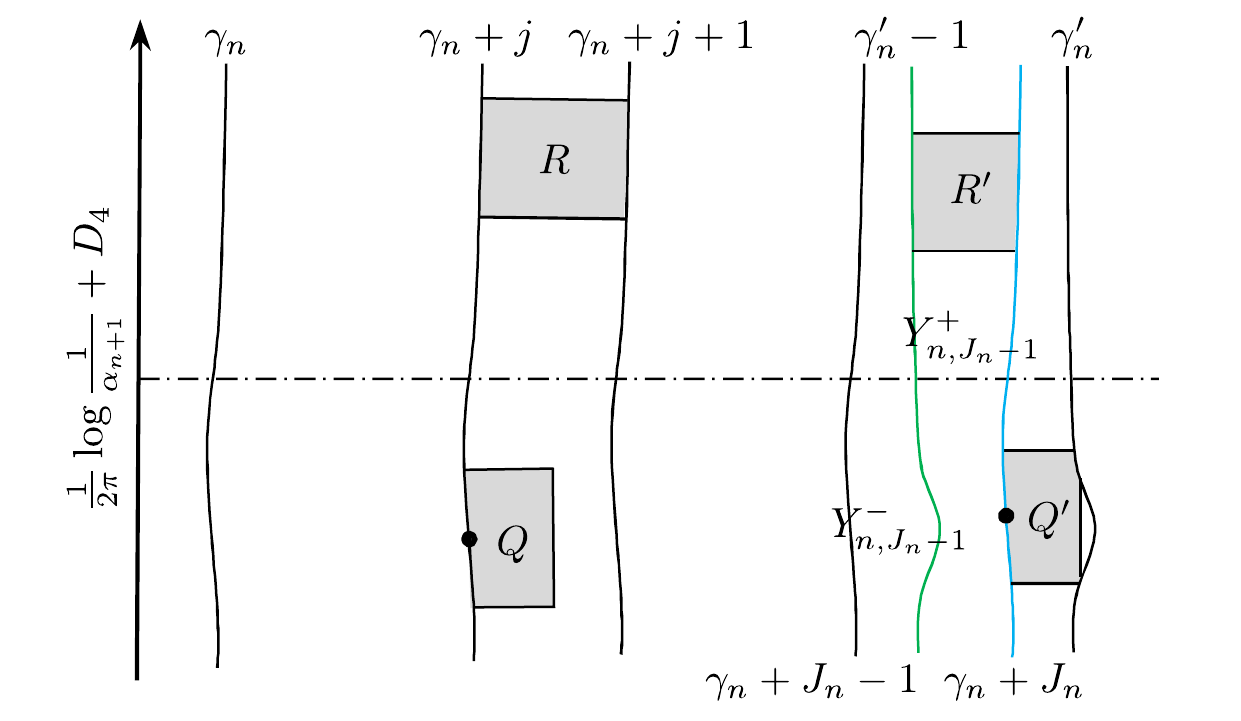}
  \caption{An sketch of two generic almost rectangles $R, R'$, and two generic nice half boxes $Q, Q'$. 
  Here $R, Q \subset Y_{n,j}$, $R' \subset Y_{n,J_n-1}^+$, and $Q' \subset Y_{n,*}$.}
  \label{Fig:almost-rectg}
\end{figure}

For $n\geq 0$, $Q \subset T_n$ is called a \textbf{nice half box} if there is $r \in [\nu_0, 3/2]$ such that either 
\begin{itemize}
\item $Q=\Boxx(\zeta,r)\cap Y_{n,j}$, for some $\zeta\in\gamma_n+j'$ with $j'\in\widetilde{\J}_n$ 
and $j \in \J_n\cap\{j'-1,j'\}$, or 
\item $Q=\Boxx(\zeta,r)\cap Y_{n,*}$, for some $\zeta \in \gamma_n' \cup (\gamma_n+J_n)$, or 
\item $Q=\Boxx(\zeta,r)\cap Y_{n,J_n-1}^-$, for some $\zeta \in (\gamma_n+J_n-1) \cup (\gamma_n'-1)$, or 
\item $Q=\Boxx(\zeta,r)\cap Y_{n,J_n-1}^+$, for some $\zeta \in (\gamma_n+J_n) \cup (\gamma_n'-1)$.
\end{itemize}
However, when $\im \zeta \leq \tfrac{1}{2\pi} \log \tfrac{1}{\alpha_{n+1}}+D_4$, we take $r=\nu_0$ in the above cases. 
That is because, since $Q \subset T_n$ and the lower parts of $T_n$ are ``bands'' of width $\nu_0$, 
there is no gain in considering larger values of $r$. 
We do this for the simplicity of later arguments. 
On the other hand, for $\zeta \in T_n$ with large enough imaginary part, we allow larger values of $r$ up to $3/2$.
In particular, some nice half boxes are also almost rectangles. 
See Figure \ref{Fig:almost-rectg} for an illustration of almost rectangles and nice half boxes. 
%Note that a nice half box may also be an almost rectangle.  

%If $\varepsilon_{n+1}=+1$, then $\Boxx(\zeta,\nu_0)\cap (\gamma_n+\N)$ with $\zeta\in\gamma_n'$ may be non-empty ($\nu_0$ is small but the width of $Y_{n,*}$ might be smaller). We will consider the images of the above two kinds of topological disks (almost rectangles and nice half boxes) under $\xi_n$ and use these two types of topological disks to pack the images.

\subsection{Good packings}\label{subsec:dens-I} 
Let $\Omega \subset \overline{Y}_n$, $n\geq 0$, be a bounded measurable set. 
By a \textbf{packing of} $\Omega$ in $T_n$, we mean any finite collection  
\[\Pack(\Omega)=\{V_{n,i}:1\leq i\leq b_n\},\]
where each $V_{n,i}$ is either an almost rectangle or a nice half box in $\Omega \cap T_n$, 
and $\area(V_{n,i}\cap V_{n,j})=0$ whenever $i \neq j$. 
Evidently, for a given $\Omega$, there may be no such collection or there may be many such collections.
%Slightly abusing notations for the sake of simplicity, let us write 
%$\area(\Pack(\Omega))=\sum_{i=1}^{b_n}\area(V_{n,i})$. 

%Recall that the \textit{density} of $\Pack(\Omega)$ in $\Omega$ is defined as 
%\begin{equation}
%\dens(\Pack(\Omega),\Omega)=\frac{\area(\Pack(\Omega))}{\area(\Omega)}.
%\end{equation}

As in Section~\ref{sec-criterion}, we use $\widehat{\Pack}(\xi_n(S))$ to denote 
the union of all sets in $\Pack(\xi_n(S))$. 

\begin{lemma}\label{lemma:image}
There is $\widetilde{\delta}>0$ such that for any $\varepsilon \in (0,1/10)$ and any $n\geq 1$ the following holds. 
For any almost rectangle or nice half box $S \subset T_{n-1}=T_{n-1}(\varepsilon)$, there exists 
$\Pack(\xi_n(S))=\{V_{n,i}:1\leq i\leq b_n\}$ of $\xi_n(S)$ in $T_n=T_n(\varepsilon)$ satisfying the following 
properties: 
\begin{itemize}
\item[(a)] $\dens(\widehat{\Pack}(\xi_n(S)), \xi_n(S))\geq\widetilde{\delta}$, 
\item[(b)] if $S$ is a nice half box, but not an almost rectangle, of height $2r \in [2\nu_0,3]$, and for all 
$\zeta \in S$, $\im \zeta \geq \tfrac{1}{2\pi}\log\tfrac{1}{\alpha_n}+D_4$ and 
$\im\xi_n(\zeta)\geq \tfrac{1}{2\pi}\log\tfrac{1}{\alpha_{n+1}}+D_4$, then every $V_{n,i}$ in $\Pack(\xi_n(S))$ 
is either an almost rectangle or a nice half box of height $\min\{8r/3,3\}$,  
\item[(c)] if $S$ is an almost rectangle such that for all $\zeta\in S$, 
$\im\zeta\geq \tfrac{1}{2\pi}\log\tfrac{1}{\alpha_n}+D_4$ and 
$\im\xi_n(\zeta)\geq \tfrac{1}{2\pi}\log\tfrac{1}{\alpha_{n+1}}+D_4$, then every $V_{n,i}$ in $\Pack(\xi_n(S))$ 
is an almost rectangle, and 
\[\dens(\widehat{\Pack}(\xi_n(S)), \xi_n(S))\geq 1-\varepsilon/5.\]
\end{itemize}
\end{lemma}

Note that each $V_{n,i}$ in $\Pack(\xi_n(S))$ obtained from Lemma~\ref{lemma:image} is contained in $T_n$. 
In particular, if some $V_{n,i}$ crosses the half-plane $\im\zeta < \tfrac{1}{2\pi}\log\tfrac{1}{\alpha_{n+1}}+D_4$, 
then $V_{n,i}$ must have height $2\nu_0$. 
Evidently, $\Pack(\xi_n(S))$ is not unique, but any such choice works for our purpose. 

We shall use the classic Koebe distortion theorem, %(see \cite[Theorem 1.6]{Pom75}), 
which states that for any univalent map $f:B_1(0) \to \C$ satisfying $f(0)=0$ and $f'(0)=1$, for every $z \in B_1(0)$, 
\[|z|/(1+|z|)^2 \leq |f(z)| \leq |z|/(1-|z|)^2, \qquad (1-|z|)/(1+|z|)^3 \leq |f'(z)| \leq (1+|z|)/(1-|z|)^3.\]

\begin{proof}
Let us present the proof when $\varepsilon_n=-1$ and $\varepsilon_{n+1}=+1$, as it simplifies 
some indexes. The other cases can be handled by minor modifications of the argument.

\medskip 

{\em (a)}
The proof is divided into several cases, depending on the locations of $S$ and $\xi_n(S)$. 

\smallskip 

{\em Case 1:} $S=\Boxx(\zeta_{n-1},\nu_0) \cap Y_{n-1,j}$ is a nice half box, 
with $\zeta_{n-1} \in \gamma_{n-1} + j$ and $j \in \J_{n-1}$. 
By periodicity, we may assume that $j=0$. By Lemma~\ref{lemma:fixed-width}--(c), 
$\xi_n:S\to Y_n$ has a univalent extension onto $\Boxx(\zeta_{n-1},20\nu_0)$.
Then, by Lemma~\ref{lemma:diameter-n} and the Koebe distortion theorem, for any 
$\zeta \in \partial \Boxx(\zeta_{n-1}, \nu_0)$, 
\begin{equation}\label{equ:est-space}
\big|\widetilde{\xi}_{n,0}(\zeta)-\widetilde{\xi}_{n,0}(\zeta_{n-1})\big|\geq \frac{\nu_0}{(1+\sqrt{2}/20)^2}\cdot\frac{5}{3}
>1.45\,\nu_0>\sqrt{2}\,\nu_0. 
\end{equation}
Therefore, $\xi_n(S)$ contains $\Boxx(\zeta_n,\nu_0)\cap Y_{n,0}$ with $\zeta_n=\xi_{n,0}(\zeta_{n-1})\in\gamma_n$. 
By \eqref{equ:T-n-above}, $\xi_n(S) \subset \BH_{h_n}$. 
By the Koebe distortion theorem, $\xi_n(S)$ has bounded geometry and hence there are 
$\widetilde{\delta}_1>0$ depending only on $\nu_0$, and $\Pack(\xi_n(S))$ such that 
$\dens(\widehat{\Pack}(\xi_n(S)), \xi_n(S))\geq\widetilde{\delta}_1$.

The same argument applies if $S=\Boxx(\zeta_{n-1},\nu_0)\cap Y_{n-1,*}$ is a nice half box for some 
$\zeta_{n-1}\in\gamma_{n-1}+J_{n-1}$, or $S=\Boxx(\zeta_{n-1},\nu_0)\cap (Y_{n-1,*}-1)$ is a nice half box for some 
$\zeta_{n-1}\in\gamma_{n-1}+J_{n-1}-1$, because $\xi_n:S\to Y_n$ has a univalent extension 
$\widetilde{\xi}_{n,*}:\Boxx(\zeta_{n-1},20\nu_0)\to\Pi_n$. 

\smallskip 

{\em Case 2:} $S=\Boxx(\zeta_{n-1}, \nu_0)\cap Y_{n-1,j}$ is a nice half box, where $\zeta_{n-1} \in \gamma_{n-1}+j+1$ 
and $j \in \J_{n-1}$. 
Without loss of generality, we assume that $j=0$. 
Here $\xi_n:S\to Y_n$ has a univalent extension onto $\Boxx(\zeta_{n-1},20\nu_0)$  
and we may repeat \eqref{equ:est-space}, with $\zeta_n=\xi_{n,0}(\zeta_{n-1})\in\gamma_n'$. 
Let us choose $\zeta_n'\in\gamma_n+J_n$ with $\im\zeta_n'=\im\zeta_n$. 
By \eqref{equ:small-arg} and \eqref{equ:est-space}, there is $\varrho>0$ and $r\geq \nu_0$ such that 
\begin{itemize} 
\item if $\re(\zeta_n-\zeta_n')\geq\varrho$, then $\xi_n(S)$ contains at least one nice half box 
$\Boxx(\zeta_n,\nu_0) \cap Y_{n,*}$, %with $r\geq\nu_0$;
\item if $\re(\zeta_n-\zeta_n')<\varrho$, then $\xi_n(S)$ contains at least one nice half box $\Boxx(\zeta_n',\nu_0) 
\cap Y_{n,J_n}$. % with $r\geq\nu_0$.
\end{itemize}
In both cases, by the Koebe distortion theorem, $\xi_n(S)$ has bounded geometry and hence there is 
$\widetilde{\delta}_2>0$ and a $\Pack(\xi_n(S))$ in $T_n$ satisfying  
$\dens(\widehat{\Pack}(\xi_n(S)), \xi_n(S))\geq\widetilde{\delta}_2$.

The same argument applies if $S=\Boxx(\zeta_{n-1},\nu_0)\cap Y_{n,J_n-1}^+$ is a nice half box with
$\zeta_{n-1}\in\gamma_{n-1}+J_{n-1}$. 

\smallskip

{\em Case 3:} $S=\Boxx(\zeta_{n-1}, \nu_0)\cap (Y_{n-1,*}-j)$ is a nice half box with 
$\zeta_{n-1} \in (\gamma_{n-1}'- j)$ with $j \in \{0,1\}$. 
First assume that $\zeta_{n-1} \in \gamma_{n-1}'$. 
Then $\xi_n:S\to Y_n$ has a univalent extension to $\widetilde{\xi}_{n,*}:\Boxx(\zeta_{n-1},20\nu_0)\to\Pi_n$ 
and \eqref{equ:est-space} holds, where 
$\zeta_n=\widetilde{\xi}_{n,*}(\zeta_{n-1})\in\gamma_n'-1$ (see Lemma~\ref{lemma:tiled}(a)). Therefore, $\xi_n(S)$
contains at least one nice half box $\Boxx(\zeta_n,\nu_0) \cap Y_{n,J_n-1}^-$. 
By the Koebe distortion theorem, $\xi_n(S)$ has bounded geometry, and hence, there is $\widetilde{\delta}_3>0$ 
and a $\Pack(\xi_n(S))$ in $T_n$ satisfying $\dens(\widehat{\Pack}(\xi_n(S)), \xi_n(S)) \geq \widetilde{\delta}_3$.

A similar argument applies when $\zeta_{n-1} \in \gamma_{n-1}'-1$. 

%By making the constants $\widetilde{\delta}_1, \widetilde{\delta}_2, \widetilde{\delta}_3$ smaller 
%if necessary, we may assume that the lower bounds on the densities hold also for the nice half boxes 
%$S=\Boxx(\zeta_{n-1},\nu_0)\cap A_{n-1}$, with $A_{n-1}$ any of the sets $Y_{n-1,J_{n-1}}^\pm$ or $Y_{n-1,j}$ 
%for $j \in \J_{n-1} \cup \{*\}$, and $S$ lies in the half plane $\im \zeta \geq \frac{1}{2\pi}\log\frac{1}{\alpha_n}+D_4$. 
%For that one employs Lemma \ref{lemma:fixed-width}-(d). 

\smallskip 

{\em Case 4:} $S \subset Y_{n-1,j}$ is an almost rectangle, with $j \in \J_{n-1}$. 
By Lemma~\ref{lemma:fixed-width}-(d) and the Koebe distortion theorem, there are $\widetilde{\delta}_4>0$ 
and $\Pack(\xi_n(S))$, with elements in $T_n$, satisfying $\dens(\Pack(\xi_n(S)), \xi_n(S))\geq\widetilde{\delta}_4$.

Similarly, the conclusion also holds if $S$ is an almost rectangle contained in $Y_{n-1,*}$, $Y_{n-1,J_{n-1}}^-$
or $Y_{n-1,J_{n-1}}^+$, or $S$ is a nice half box in the half plane 
$\im \zeta \geq \tfrac{1}{2\pi} \log \tfrac{1}{\alpha+{n+1}}+ D_4$. 

We may define $\widetilde{\delta}=\min\{\widetilde{\delta}_i:1\leq i\leq 4\}$.

\medskip

{\em (b)} Given $S= \Boxx(\zeta_{n-1},r)$ as in (b), 
%\cap \star_{n-1}$, where $r \in [\nu_0,3/2]$ and $\star_{n-1}$ is any of $Y_{n-1,J_{n-1}-1}^\pm$ or $Y_{n-1,j}$ with $j \in \J_{n-1} \cup \{*\}$. Suppose that $\im\zeta\geq \tfrac{1}{2\pi}\log\tfrac{1}{\alpha_n}+D_4$ and $\im\xi_n(\zeta)\geq \tfrac{1}{2\pi}\log\tfrac{1}{\alpha_{n+1}}+D_4$ for all $\zeta\in S$. 
by Lemma~\ref{lemma:diameter-n}, there is $\Pack(\xi_n(S))$ whose elements are almost rectangles or nice half boxes of the form 
$\Boxx(\zeta_n,\min\{4r/3,3/2\})\cap \star_n$, where $\zeta_n\in\gamma_n+\widetilde{\J}_n\cup\{*,*-1\}$ and `$\star_n$' 
denotes $Y_{n,J_n-1}^\pm$ or $Y_{n,j}$ with $j\in\J_n\cup\{*\}$. 

\smallskip

{\em (c)} Given $S$ as in (c), $\xi_n:S \to \Pi_n$ has a univalent extension 
$\widetilde{\xi}_n:B_{\nu_0}(S)\to \Pi_n$ by Lemma \ref{lemma:fixed-width}-(d). 
By Lemma \ref{lemma:chi-inverse}-(a), for any $\zeta_{n-1}\in S$ and $\zeta_n=\xi_n(\zeta_{n-1})$, 
%\begin{equation}
$|\arg\widetilde{\xi}_n'(\zeta_{n-1})| \leq \varepsilon/10$.
%\end{equation}
Because the height of each almost rectangle is at least $+1$, $\xi_n(S)$ can be packed by a family of 
almost rectangles $\{V_{n,i}:1\leq i\leq b_n\}$ in $T_n$ such that $\dens(\widehat{\Pack}(\xi_n(S)), \xi_n(S))\geq 1-\varepsilon/5$.
\end{proof}

%\begin{rmk}
%The proof of Lemma \ref{lemma:image} can be simplified a lot if we assume that $N$ is larger than a very big integer. In that case,
%the image of any nice half box $Q_{k-1}$ will cover at least one almost rectangle if $Q_{k-1}$ is above the height
%$\tfrac{1}{2\pi}\log\tfrac{1}{\alpha_k}+D_4$ and $\xi_k(Q_{k-1})$ is also above the height
%$\tfrac{1}{2\pi}\log\tfrac{1}{\alpha_{k+1}}+D_4$, where $k\geq 1$.
%\end{rmk}

\subsection{The nest} \label{subsec:nest}
%The main strategy to build the nest of sets is to first successively apply Lemma~\ref{lemma:image} to 
%define a sequence of collections $\{\MF_n\}_{n\geq 0}$ such that each $\MF_n= \{F_{n,i}, 1 \leq i \leq l_n \}$ consists of 
%almost rectangles and nice half boxes in $T_n$. Then we use the changes of coordinates $\{\chi_j\}_{j\geq 0}$ 
%to map those elements $F_{n,i}$ to the top level to define a nest of sets $\{\MK_n\}_{n\geq 0}$. 
%The details are presented below. 

Let us define 
\[K_{0,1} = F_{0,1}=\{\zeta \in Y_{0,0} \mid h_0 \leq \im \zeta \leq h_0+1\}, \qquad 
\MK_0  = \MF_0=\{K_{0,1}\}.\]
Then, $F_{0,1}$ is an almost rectangle. 
By Lemma~\ref{lemma:image}, there is $\MF_1:=\Pack(\xi_1(F_{0,1}))=\{F_{1,i}:1\leq i\leq l_1\}$ which 
is not empty and consists of almost rectangles and half boxes in $T_1$. 
Applying Lemma~\ref{lemma:image} to each $F_{1,i}$, for $1\leq i \leq l_1$, we obtain 
\[\MF_{2}^{1,i}:= \Pack(\xi_2(F_{1,i}))= \{F_{2,k}^{1,i}, 1 \leq k \leq l_{2}^{1,i}\}.\]
Combining these collections, we obtain $\MF_2= \cup_{1 \leq i \leq l_1}\MF_2^{1,i}$, 
which consists of almost rectangles and half boxes. 
Note that the elements of $\MF_2$ might have overlaps. 

Successive applications of Lemma~\ref{lemma:image} in the same fashion produce the collections 
\[\MF_n^{n-1,i}=\Pack(\xi_n(F_{n-1,i}))=\{F_{n,k}^{n-1,i}:1\leq k\leq l_n^{n-1,i}\}, \quad 1 \leq i \leq l_{n-1} \]
which are combined to define 
\[\MF_n= \{F_{n,k}^{n-1,i}:1\leq i\leq l_{n-1},1\leq k\leq l_n^{n-1,i}\} = \{F_{n,j}:1\leq j\leq l_n\}, \quad \forall n \geq 0\]
%This completes the definition of the collections $\MF_n$, for $n\geq 0$. 
Each $F_{n,j}$ is either an almost rectangle or a nice half box contained in $T_n$. 

For $n=1$, we let 
\[\MK_n= \MK_1=\{K_{1,i}=\xi_1^{-1}(F_{1,i})=\chi_{1,(\varepsilon_1+1)/2}(F_{1,i}):1\leq i\leq l_1\}.\]
Fix $n \geq 2$. 
For each $F_{n,i_n}$, with $1 \leq i_n \leq l_n$, there exists a unique 
sequence $(i_0,i_1,\cdots,i_{n-1})$, satisfying $1\leq i_m\leq l_m$ and $0\leq m\leq n-1$, such that
\[F_{m+1,i_{m+1}}\in\Pack(\xi_{m+1}(F_{m,i_m})).\]
Then, for each $m$ satisfying $0 \leq m \leq n-1$, there is an inverse branch of $\xi_{m+1}$ that maps 
$F_{m+1,i_{m+1}}$ to $F_{m,i_m}$.
Using those inverse branches, we define
\[\MK_n= \{\xi_1^{-1}\circ\cdots\circ\xi_n^{-1}(F_{n,i_n}): 1\leq i_n \leq l_n\} = \{K_{n,j}:1\leq j\leq l_n\}.\]
We note that $l_n= \textstyle{\sum_{i=1}^{l_{n-1}}} l_n^{n-1,i}$.

It follows from the construction that 
\begin{itemize}
\item $\area(K_{n,i} \cap K_{n,j})=0$, for all $1 \leq i < j \leq l_n$, and 
\item for all $n\geq 1$, $\bigcup_{1 \leq j \leq l_{n+1}} K_{n+1,j} \subset  \bigcup_{1 \leq j \leq n} K_{n,j}$.  
\end{itemize}

\section{Dimension of the post-critical set}\label{sec:H-dim-PC}
In this section we prove Theorem~\ref{thm:dim-pc}, which is a generalisation of 
Theorem~\ref{thm:quadratic-dim-pc}.
The proof is based on the criterion presented in Section~\ref{sec-criterion}, which we apply to the nest 
$\{\MK_n\}_{n\geq0}$ built in Section~\ref{sec-partition}. 
In Section~\ref{sec:densities-diameters} we estimate the relative densities and diameters for the nest, and in 
Section~\ref{subsec:proof-dim-2} we combine those estimates to present the proof. 

\subsection{Densities and diameters}\label{sec:densities-diameters}
Given a univalent or anti-univalent map $g: \Omega \to \mathbb{C}$, where $\Omega \subset \C$, 
we define the \textbf{distortion} of $g$ on $\Omega$ as 
\[L(g|_\Omega)= \sup \{|g'(x)/g'(y)| : x,y \in \Omega\}. \] 
We say that $g$ has \textbf{bounded distortion} on $\Omega$ if $L(g|_\Omega)$ is finite. 
Evidently, for any pair of univalent or anti-univalent maps $g_1:\Omega_1\to \C$ and $g_2:\Omega_2\to\C$ 
with $\Omega_1$ and $\Omega_2$ bounded and $g_1(\Omega_1)\subset\Omega_2$, we have 
\begin{equation}\label{equ:distor-1}
L(g_1|_{\Omega_1})=L(g_1^{-1}|_{g_1(\Omega_1)}), \quad 
L((g_2\circ g_1)|_{\Omega_1})\leq L(g_1|_{\Omega_1})L(g_2|_{g_1(\Omega_1)}).
\end{equation} 
Also, for any measurable set $X \subset \Omega$, 
\begin{equation}\label{equ:density-inequ}
L(g|_\Omega)^{-2}\,\dens(g(X),g(\Omega))\leq \dens(X,\Omega)\leq L(g|_\Omega)^2\,\dens(g(X),g(\Omega)).
\end{equation}

\begin{lemma}\label{lemma:distortion}
There is a constant $M_3$ such that for all $n\geq 1$ and $1 \leq i \leq l_{n-1}$, 
\[ L\left (\xi_n\circ\cdots\circ\xi_1:K_{n-1,i}\to \xi_n(F_{n-1,i}) \right ) \leq M_3.\]
\end{lemma}

\begin{proof}
For $1\leq i\leq l_{n-1}$, each $F_{n-1,i}$ is an almost rectangle or a nice half box. 
By Lemma~\ref{lemma:fixed-width}--(c),(d), and \eqref{equ:small-arg}, 
$\xi_n:F_{n-1,i}\to \xi_n(F_{n-1,i})$ has a univalent or anti-univalent extension 
onto $B_{\nu_0}(F_{n-1,i})$.
Because each nice half box and almost rectangle has a uniformly bounded diameter, 
the conformal modulus of $B_{\nu_0}(F_{n-1,i})\setminus F_{n-1, i}$ is uniformly bounded from below.
Then, by the Koebe distortion theorem, there exists $M_3'$, independent of $n$ and $i$, 
such that $L(\xi_n|_{F_{n-1,i}})\leq M_3'$.

On the other hand, $G_{n-1}^{-1}=(\xi_{n-1} \circ \cdots\circ\xi_1)^{-1}:F_{n-1,i}\to K_{n-1,i}$ has a univalent or 
anti-univalent extension to $\widetilde{G}_{n-1}^{-1}:B_{\nu_0}(F_{n-1,i})\to \C$.
%Let $\widetilde{K}_{n-1,i}=\widetilde{G}_{n-1}^{-1}(B_{\nu_0}(F_{n-1,i}))$.
%Then $K_{n-1,i}\subset \widetilde{K}_{n-1,i}$ and
%\begin{equation}
%\Mod(\widetilde{K}_{n-1,i}\setminus K_{n-1,i})\geq \kappa.
%\end{equation}
Applying the Koebe distortion theorem in the same fashion, there must be $M_3''$, independent of $n$ and $i$, 
such that $L(G_{n-1}^{-1}|_{F_{n-1,i}})\leq M_3''$.
The desired bound in the lemma follows from \eqref{equ:distor-1}, by setting $M_3=M_3'M_3''$.
\end{proof}

%\begin{rmk}
%In \cite[Section 3]{McM87}, McMullen controlled the distortion by studying the size of the boxes, the amount of expansion and the nonlinearities.
%\end{rmk}

Recall $D_4=D_4(\varepsilon)$ from Lemma \ref{lemma:chi-inverse} and $\widetilde{\mu}_n$ from 
Lemma~\ref{lemma:diameter-n}. Also recall that $\hat{\MK}_n= \cup_{K \in \MK_n} K$. 

\begin{cor}\label{cor:density-low}
There are $\delta \in (0,1)$ and $M_4 \geq 1$ such that for any $\varepsilon \in(0,1/10)$, the following hold: 
\begin{enumerate}
\item For all $n\geq 1$ and all $1\leq i\leq l_{n-1}$,
\[\dens(\hat{\MK}_n, K_{n-1,i})\geq\delta.\]
\item If $F_{n-1, i}$ is an almost rectangle and for all $\zeta \in F_{n-1,i}$,
\[\im \zeta \geq \tfrac{1}{2\pi}\log\tfrac{1}{\alpha_n}+D_4, \quad \im \xi_n(\zeta) \geq \tfrac{1}{2\pi}\log\tfrac{1}{\alpha_{n+1}}+D_4,\] 
then 
\[\dens(\hat{\MK}_n, K_{n-1,i})\geq 1-M_4\varepsilon.\]
\item For all $n\geq 1$ and all $1\leq i\leq l_n$, 
\[\diam(K_{n,i})\leq M_4 \textstyle{\prod_{k=1}^n} \widetilde{\mu}_k.\]
\end{enumerate}
\end{cor}

\begin{proof}
Part (a) follows from Lemmas \ref{lemma:image} and \ref{lemma:distortion}, as well as \eqref{equ:density-inequ}. 
Similarly, for (b), we note
\[\dens(\hat{\MK}_n, K_{n-1,i})\geq 1-\dens(\C\setminus \hat{\MK}_n, K_{n-1,i})\geq 1-M_3^2\varepsilon/5.\]
Every $F_{n,i}=\xi_n\circ\cdots\circ\xi_1(K_{n,i})\subset T_n$, for $n\geq 1$ and $1 \leq i \leq l_n$, 
is an almost rectangle or a nice half box with uniformly bounded diameter. 
Then, (c) follows from Lemmas~\ref{lemma:diameter-n} and \ref{lemma:distortion}.
\end{proof}

\subsection{Main theorem}\label{subsec:proof-dim-2}
Recall the class of maps $\QIS_\alpha$ from Section~\ref{subsec-IS-class}, and recall that $N$ is determined in \eqref{E:N-condition}. 

\begin{thm}\label{thm:dim-pc}
For every $\alpha \in \HT_N \setminus \HH$ and every $f \in \QIS_\alpha$, $\dim_H (\Lambda(f))=2$.
\end{thm}

\begin{proof}
The class of maps $\QIS_\alpha$ is defined in Section~\ref{subsec-IS-class}, and $N$ is determined
in Section~\ref{subsec:renorm-tower}. 
Recall $M_4>1$ and $\delta \in (0,1)$ from Corollary \ref{cor:density-low}. 
Fix an arbitrary $\varepsilon$ satisfying $0< \varepsilon<(1-\delta)/(10M_4)$. 
Let $D_4 = D_4(\varepsilon)$ and $D_4'=D_4'(\varepsilon)$ be the constants obtained in Lemma \ref{lemma:chi-inverse}, 
and let $h_n=(4/3)^n D_4'$, for $n\geq 0$, as in \eqref{equ:h-n}.
We shall use the nest $(\MK_n)_{n\geq 0}$ and the collections $(\MF_n)_{n\geq 0}$, 
defined in Section~\ref{subsec:nest}. 
Recall that all elements of $\MF_n$ are in $T_n \subset \mathbb{H}_{h_n}$. 

%We will compare $h_n$ with $\tfrac{1}{2\pi}\log\tfrac{1}{\alpha_{n+1}}+D_4$ and divide the arguments into several cases.

By Lemma~\ref{lemma:image}-(b), there is $I \geq 1$ such that for any $n \geq 0$, if for all $1 \leq j \leq I$, we have 
\[h_{n+j-1} \geq \tfrac{1}{2\pi}\log\tfrac{1}{\alpha_{n+j}}+D_4, \]
then any element of $\MF_{n+I-1}= \{F_{n+I-1,i}:1 \leq i\leq l_{n+I-1}\}$ is an almost rectangle.
In other words, if we start with a nice half box (or an almost rectangle) on some level $n$, which remain above those 
heights as we map deeper in the renormalisation tower, the sizes of the nice half boxes in the packing increase, 
and after at most $I$ levels, they all become almost rectangle. 
We need to consider this because, once we reach almost rectangles, we may use the better estimate in 
Corollary~\ref{cor:density-low}-(b). 

Recall $M_2$ from Lemma~\ref{lemma:diameter-n}. 
We define $(\mu_k)_{k\geq 1}$ and $(\delta_k)_{k\geq 1}$ according to the following: 
\begin{itemize}
\item[(L1)] If $h_{k-1}< \tfrac{1}{2\pi}\log\tfrac{1}{\alpha_k}+D_4$, let $\mu_k=M_2 e^{-2\pi h_{k-1}}$ 
and $\delta_k=\delta$.
\item[(L2)] If $h_{k-1} \geq \tfrac{1}{2\pi}\log\tfrac{1}{\alpha_k}+D_4$, let $m_1 \geq 1$ be the smallest integer 
and $m_2 \geq  k$ be the largest integer such that for all $j$ satisfying $m_1 \leq j \leq m_2$, we have 
$h_{j-1} \geq \tfrac{1}{2\pi}\log\tfrac{1}{\alpha_j}+D_4$. 
Then, for  $m_1 \leq j \leq m_2$, let 
\[\mu_j=3/5, \qquad \delta_j= \begin{cases}
\delta & \text{ if } m_1 \leq j \leq m_1 + I-2 , \text{ or } j=m_2, \\
1- M_4 \varepsilon & \text{ if }  m_1+ I-1 \leq j \leq m_2-1.  
\end{cases}   \]
\begin{comment}
\item[(L2)]  If $h_{k+j-1}\geq\tfrac{1}{2\pi}\log\tfrac{1}{\alpha_{k+j}}+D_4$ for $0\leq j\leq m$ 
with $0\leq m\leq I-1$, $h_{k-2}<\tfrac{1}{2\pi}\log\tfrac{1}{\alpha_{k-1}}+D_4$ and 
$h_{k+m}<\tfrac{1}{2\pi}\log\tfrac{1}{\alpha_{k+m+1}}+D_4$, we define 
\[\mu_{k+j}=\frac{3}{5} \text{\quad and\quad} \delta_{k+j}=\delta, \text{\quad where } 0\leq j\leq m.\]
\item[(L3)] If $h_{k+j-1}\geq\tfrac{1}{2\pi}\log\tfrac{1}{\alpha_{k+j}}+D_4$ for $0\leq j\leq m$ 
with\footnote{Actually, $m$ cannot be $+\infty$ if $\alpha$ is not of Herman type.} $I\leq m\leq +\infty$ and
$h_{k-2}<\tfrac{1}{2\pi}\log\tfrac{1}{\alpha_{k-1}}+D_4$, we define
\[\mu_{k+j}=\frac{3}{5} \text{ for } 0\leq j\leq m,
\text{\quad and\quad}
\delta_{k+j}=
\begin{cases} 
\delta   & \text{ if } 0 \leq j \leq I-1,\\
1-M_4 \varepsilon & \text{ if } I\leq j\leq m. 
\end{cases}\]
\end{comment}
\end{itemize} 

Since $\alpha$ is not a Herman number, one can show that, in (L2), $m_2$ must be finite (See for instance \cite{Che13}). 
However, we do not need the finiteness of $m_2$ here (see the final paragraph in this proof). 
 
By our choice of $\varepsilon$, we have $\delta \leq 1-M_4 \varepsilon$. 
Then, it follows from Lemma~\ref{lemma:image} and Corollary~\ref{cor:density-low} that for all $n \geq 1$ and all $1 \leq i \leq l_{n-1}$,
\[\dens(\hat{\MK}_n, K_{n-1,i})\geq\delta_n, \qquad \diam(K_{n,i}) \leq d_n:=M_4\textstyle{\prod_{k=1}^n} \mu_k.\]
We claim that
\[\limsup_{n\to\infty} \frac{\sum_{k=1}^{n+1}|\log \delta_k |}{\sum_{k=1}^n|\log \mu_k| -\log M_4} \leq 4M_4\varepsilon.\]
Because $\sum_{k=1}^\infty |\log \delta_k|=\sum_{k=1}^\infty |\log \mu_k| = +\infty$, 
and $| \log \delta_k| \in [0, \log (1/\delta)]$, 
%$|\log\delta_k|\in [0,\log(1/\delta)]$, for all $k\geq 1$, and $0< \delta<1-M_4\varepsilon$, 
it suffices to show that 
\[\limsup_{n\to\infty} \frac{\sum_{k=1}^n|\log \delta_k |}{\sum_{k=1}^n|\log \mu_k|} \leq 4M_4\varepsilon.\]
We deal with that by considering two cases. 

\textbf{(i)} There are only finitely many distinct $k$ satisfying $L1$. 
%integers $1\leq k_1<k_2<\cdots<k_\ell$ such that  $h_{k_i-1}<\tfrac{1}{2\pi}\log\tfrac{1}{\alpha_{k_i}}+D_4$, for $1\leq i\leq \ell$. This means that 
%\begin{equation}
%\begin{split}
%\big|\log\mu_{k_i}\big|=\log\frac{e^{2\pi h_{k_i-1}}}{M_2} &\text{\quad for\quad} 1\leq i\leq \ell\text{\quad and\quad} \\
%\big|\log\mu_k\big|=\log\frac{5}{3} &\text{\quad for\quad} k\not\in\{k_i:1\leq i\leq\ell\}.
%\end{split}
%\end{equation}
Then, for large enough $k$, we have $\mu_k=3/5$ and $\delta_k=1- M_4 \varepsilon$. 
Therefore, using $M_4 \varepsilon <1/10$ and $|\log (5/3)| \geq 1/2$, 
%$k\geq k_\ell+I$, we have $\log\delta_k=-\log(1-M_4\varepsilon)$. This implies that
\[\lim_{n\to\infty} \frac{\sum_{k=1}^n|\log \delta_k |}{\sum_{k=1}^n|\log \mu_k|}  
= \frac{\log(1-M_4\varepsilon)}{\log(3/5)} \leq (10/9) M_4 \varepsilon \times 2 < 4M_4\varepsilon. \]

\textbf{(ii)} There are infinitely many distinct integers $k_1<k_2< k_3  < \cdots$ such that
$h_{k_i-1}<\tfrac{1}{2\pi}\log\tfrac{1}{\alpha_{k_i}}+D_4$, and $h_{k-1}\geq \tfrac{1}{2\pi}\log\tfrac{1}{\alpha_k}+D_4$ for
$k\not\in\{k_i:i\geq 1\}$. 
%Then, $\left | \log \mu_k \right | = \log (e^{2\pi h_{k-1}}/M_2)$ if $k \in \{k_i : i\geq 1\}$, 
%and $\left | \log\mu_k\right| = \log(5/3)$ if $k\not\in\{k_i:i\geq 1\}$.

To simplify the presentation, let $k_0=0$ (without assuming that $k_0$ is among the list of $k_i$). 
For any $j\geq 1$, 
\[u_j:= \textstyle{\sum_{i=k_{j-1}+1}^{k_j}} |\log \delta_i | \leq  I | \log \delta | + (k_j-k_{j-1}-1) | \log (1-M_4\varepsilon)|,\]
\[v_j:= \textstyle{\sum_{i=k_{j-1}+1}^{k_j}} |\log \mu_i | = (k_j-k_{j-1}-1) \log (5/3) + | \log (M_2 e^{-2\pi h_{k_j-1}}).\]
Fix an arbitrary $n \geq 1$, and let $\ell= \ell(n) \geq 1$ be such that $k_{\ell-1} \leq n < k_\ell$. In the same fashion, we have 
\[u_\ell':=\textstyle{\sum_{i=k_{\ell-1}+1}^n} |\log \delta_i | 
\leq I | \log \delta | + (n-k_{j-1}-1) |\log (1-M_4\varepsilon)|,\]
\[v_\ell':=\textstyle{\sum_{i=k_{\ell-1}+1}^n} |\log \mu_i | = (n-k_{j-1}-1) \log( 5/3 ).\]
Combining the above equations, we get 
\[\textstyle{\sum_{k=1}^n } |\log \delta_k |= \textstyle{\sum_{j=1}^{\ell-1}} u_j+ u_\ell' 
\leq \ell \,I |\log \delta| + (n-\ell)| \log (1-M_4\varepsilon)| \]
and
\[\textstyle{\sum_{k=1}^n} | \log \mu_k | = \textstyle{\sum_{j=1}^{\ell-1}} v_j+ v_\ell'
=\textstyle{\sum_{j=1}^{\ell-1}}  |\log (M_2 e^{-2\pi h_{k_j-1}})| + (n-\ell) \log (5/3).\]
Since $h_{k_j-1} \to \infty$, as $j \to \infty$, and $\ell=\ell(n) \to \infty$ as $n \to \infty$, it follows that
%\begin{equation}
%\lim_{\ell \to \infty} \frac{\ell \, I | \log \delta| } {\sum_{j=1}^{\ell-1} | \log (M_2 e^{-2\pi h_{k_j-1}})|}=0.
%\end{equation}
%Note that 
\[\begin{aligned}
\limsup_{n\to\infty} \frac{\sum_{k=1}^n|\log \delta_k |}{\sum_{k=1}^n |\log \mu_k|} 
& \leq \limsup_{n \to \infty}\frac{\ell\,I |\log \delta|}{\sum_{k=1}^n|\log \mu_k|}+
\limsup_{n\to\infty}\frac{(n-\ell) |\log(1-M_4\varepsilon)|}{\sum_{k=1}^n|\log \mu_k|} \\ 
& \leq \lim_{\ell \to \infty} \frac{\ell \, I | \log \delta| } {\sum_{j=1}^{\ell-1} | \log (M_2 e^{-2\pi h_{k_j-1}})|} 
+ \lim_{n\to\infty}\frac{(n-\ell) |\log(1-M_4\varepsilon)|}{(n-\ell) \log (5/3)} 
\leq 4M_4\varepsilon.
\end{aligned}\]

By Proposition~\ref{prop:McMullen}, we must have $\dim_H(\bigcap_{n\geq 0} \hat{\MK}_n)\geq 2-4M_4 \varepsilon$.
However, although $\varepsilon >0$ was arbitrary, we may not conclude that the nest shrinks to a set of dimension 2. 
That is because the nest itself depends on $\varepsilon$. 

The set $\Phi_0^{-1}(\bigcap_{n \geq 0} \hat{\MK}_n)$ is contained in $\Lambda (f_0) \cup \Delta(f_0)$, 
see \cite[Proposition 5.10]{Che19} or \cite[Section 7.5]{Che17}. 
Evidently, $\Phi_0^{-1}$ is conformal on a neighbourhood of $K_{0,1}$ and, by Theorem~\ref{thm-Lambda-topo}, 
when $\alpha \in \HT_N\setminus\HB$, $\Delta(f_0)=\emptyset$. 
Therefore, $\dim_H(\Lambda(f_0) \geq \dim_H(\bigcap_{n\geq 0}\hat{\MK}_n) \geq 2 - 4M_4 \varepsilon$, and hence 
$\dim_H(\Lambda(f_0) = 2$. 

Suppose that $\alpha \in \HT_N\cap(\HB\setminus\HH)$. 
Then every $f_n$, for $n \geq 0$, has a Siegel disk $\Delta(f_n)$ whose boundary does not contain the unique 
critical point of $f_n$. 
Recall the sets $\Pi_n$, $n\geq 0$, from Section~\ref{subsec:change}, and let 
\[\widetilde{\Delta}(f_n)= \{\zeta \in \Pi_n \mid \Phi_n^{-1}(\zeta) \in \Delta(f_n)\} \quad \text{ and } \quad   
y_n= \inf \{ \im \zeta \mid \zeta \in \widetilde{\Delta}(f_n)\}.\]
We claim that $\lim_{n\to\infty}y_n=+\infty$. 
Otherwise, by Corollary~\ref{cor:hyperbolic-contraction}, the critical point of $f_0$ must lie on 
$\partial \Delta(f_0)$, which contradicts Theorem~\ref{thm-Lambda-topo}. 

Now, fix an arbitrary $\varepsilon \in (0, 1/10)$, and let $D_4'$ be the constant in Lemma~\ref{lemma:chi-inverse}. 
By the above paragraph, there is $n_0 \geq 0$ such that $y_{n_0} \geq D_4'+1$. 
Now, we apply the construction to the map $f_{n_0}$, and note that 
$K_{0,1} \cap \widetilde{\Delta}(f_{n_0}) = \emptyset$. 
Therefore, $\bigcap_{n\geq 0} \hat{\MK}_n \subset \Lambda(f_{n_0})$, and hence 
$\dim_H \Lambda(f_{n_0}) \geq 2 - 4M_4 \varepsilon$. 
Through the changes of coordinates $\chi_j$, we have $\dim_H \Lambda(f_0) = \dim_H \Lambda(f_{n_0})$. 
Thus, $\dim_H \Lambda(f_0) = 2$. 
\end{proof}

Evidently, there are polynomials and rational functions with an irrationally indifferent fixed point whose 
restriction to a neighbourhoods of the fixed point belongs to $\IS$. 

\begin{cor}\label{cor:dim-2}
If $\alpha \in \HT_N \setminus \HH$ and $f$ is a rational function which belongs to the class $\IS_\alpha$, 
then the Julia set of $f$ has Hausdorff dimension two.
\end{cor}
\section{Dimension of the hairs without the end points} \label{section:HD-hairs}
In this section we prove Theorem~\ref{thm:dim-hairs} which generalises Theorems~\ref{thm:dim-hairs-quadratic-(ii)} 
and \ref{thm:dim-hairs-quadratic-(iii)}. 
Recall that $\Lambda(f)$ is the post-critical set of $f$, and $\Delta(f)$ is the Siegel disk 
of $f$ when $f$ is linearisable, and $\Delta(f)=\{0\}$ when $f$ is not linearisable. 
By Theorem~\ref{thm-Lambda-topo}, for any $f \in \QIS_\alpha$ with $\alpha\in\HT_N\setminus\HH$,  
$\Lambda(f)$ is either a Cantor bouquet or a one-sided hairy Jordan curve. 
The set $\Lambda(f) \setminus \overline{\Delta(f)}$ consists of uncountably many components, which are simple arcs. 
We also refer to each such component as a \textbf{hair} of $\Lambda(f)$. 
We say that a point $z \in \Lambda(f)$ is an \textbf{end point} of $\Lambda(f)$, if there is no continuous and injective map 
$I: (-1, +1) \to \Lambda(f)$ satisfying $I(0)=z$. 
As in the introduction, $\ME(f)$ denotes the set of all end points of $\Lambda(f)$. 
It follows that $\ME(f)$ is totally disconnected, but dense in $\Lambda(f)$. 

%The sets $\HJ$ and $\HS$ are defined in Section~\ref{section:arithmetic rot numbers}. 
Let us fix an arbitrary $f \in \QIS_\alpha$ with $\alpha \in \HT_N$.  
We shall frequently use the notations introduced in Sections \ref{subsec:renorm-tower} and \ref{subsec:change}, 
while we recall the infrequently used notations. 

\subsection{Main objects in the Fatou coordinates}\label{subsec:main-objects}
Recall the curves $\gamma_n$, $n\geq 0$, from Section~\ref{subsec:partitions-deep}. 
We need to extend each $\gamma_n$ so that it touches the boundary of $\Pi_n$. 
There exists a collection of curves $\{\tilde{\gamma}_n : [-1, 0] \to \mho\}_{n\geq 0}$ satisfying the 
following properties: 
\begin{itemize}
\item[(i)] $\tilde{\gamma}_n(-1) = 1-2i \in \partial \Pi_n$, $\tilde{\gamma}_n(0)=1= \Phi_n(-4/27)$, 
\item[(ii)] $\Phi_n^{-1}(\tilde{\gamma}_n[-1, 0]) \cap \Lambda(f_n) = \{-4/27\}$, 
\item[(iii)] $\chi_n(\tilde{\gamma}_n) \subset \tilde{\gamma}_{n-1}$. 
\end{itemize}
There are many possible choices for the above collection of curves, and any choice satisfying these properties works for our purposes.
Indeed, each such curve $\tilde{\gamma}_n$ consists of incremental pieces mapped up from deep levels 
of renormalisation tower. 
For more details one may refer to \cite{Che17} section 6.4 (where $\tilde{\gamma}_n$ is denoted by $w_n^+$). 
In particular, by the second property in the above list, $\tilde{\gamma}_n \cap \gamma_n= \{1\}$, 
and $\tilde{\gamma}_n \cup \gamma_n$ divides $\Pi_n$ into two pieces. 

As in Section~\ref{subsec:partitions-deep}, for each $\tilde{\gamma}_n$, there is a counterpart curve 
$\tilde{\gamma}_n'$ in $\Phi_n(S_n)+ \{k_n-1\}$ such that 
$\Phi_n^{-1}(\tilde{\gamma_n}')= \Phi_n^{-1}(\tilde{\gamma_n})$. 
However, in this case $\tilde{\gamma}_n'$ is not unique, but there are two choices for it due to $\Phi_n^{-1}$
having a critical point at the finite end of $\gamma_n'$. 
Out of those two choices, we choose the one which lies immediately after $\gamma_n'$ in the counter clockwise 
direction at the finite end of $\gamma_n'$. 
Then, $\gamma_n' \cup \tilde{\gamma}_n'$ divides $\Pi_n$ into two components as well. 
The reason for considering these curves is to enlarge $Y_n(D_2')$ from Section~\ref{subsec:partitions-deep} all the way down 
to the bottom boundary of $\Pi_n$. 

For $n \geq 0$, let us define our frequently used notation 
\[Y_n'\] 
to denote the closure of the connected component of $\Pi_n$ cut off by $\gamma_n\cup \tilde{\gamma}_n$ 
and $\gamma_n' \cup \tilde{\gamma}_n'$. 
It follows that $\Phi_n^{-1}$ is defined and univalent on $Y_n'$, with a unique critical point on $\partial Y_n'$. 

For $n  \geq 0$, let 
\[\Lambda_n= \left( \Phi_n^{-1} (\MP_n \cap \Lambda(f_n)) + \Z \right ) \cap Y_n'.\]
One may see that $\Lambda_n$ is the pre-image of $\Lambda(f_n)$ via the map $\Phi_n^{-1}: Y_n' \to \C$, 
although this may not be immediately clear from the definition of $Y_n'$ (and will not be used here).  
In the same fashion, for $n\geq 0$, if $f_n$ is linearisable at $0$, we define 
\[\Delta_n= \left( \Phi_n^{-1} (\MP_n \cap \Delta(f_n)) + \Z \right ) \cap Y_n',\]
and when $f_n$ is not linearisable at $0$ we let $\Delta_n=\emptyset$. 
Similarly, for $n\geq 0$, let 
\[\ME_n= \left( \Phi_n^{-1} (\MP_n \cap \ME(f_n)) + \Z \right ) \cap Y_n', \qquad 
\MH_n= \Lambda_n \setminus (\overline{\Delta}_n \cup \ME_n).\]
The set $\MH_n$ is the set of hairs without their end points on level $n$. 
It follows that $\chi_n: \Lambda_n \to \Lambda_{n-1}$, $\chi_n: \Delta_n \to \Delta_{n-1}$, 
$\chi_n: \ME_n \to \ME_{n-1}$, and $\chi_n: \MH_n \to \MH_{n-1}$. 

Recall the set $\J_n$ from Section~\ref{subsec:partitions-deep}. 
For $j \in \J_n$, let $Y_{n,j}'$ denote the closure of the component of $\Pi_n$ 
cut off by $(\gamma_n \cup \tilde{\gamma}_n)+j$ and $(\gamma_n \cup \tilde{\gamma}_n)+j+1$, and let 
$Y_{n,\diamond}$ denote the closure of the component of $Y_n'$ cut off by $\gamma_n'\cup \tilde{\gamma}_n'$ 
and $(\gamma_n' \cup \tilde{\gamma}_n')-1$. 
Using those notations, we may define $\Lambda_{n,j}= \Lambda_n \cap Y_{n,j}'$, $\Lambda_{n,\diamond}
= \Lambda_n \cap Y_{n,\diamond}'$, $\Delta_{n,j}= \Delta_n \cap Y_{n,j}'$, etc. 
In particular, we get the decomposition
\begin{equation}\label{equ:decomp-1}
\Lambda_n = \textstyle{\bigcup_{j\in\J_n\cup\{\diamond \}}} \Lambda_{n,j}.
\end{equation}

We need a uniform contraction of the maps $\chi_n$ on $Y_n'$. 
While this is guaranteed by Proposition~\ref{prop:esti-chi} at points with sufficiently large imaginary parts, 
it is not known whether it holds on all of $Y_n'$. 
To this end, we employ a uniform contraction with respect to a suitable hyperbolic metric on some 
neighbourhood of $Y_n'$, as we explain here. 
It follows from the pre-compactness of $\IS$ (see Section 7.3 in \cite{Che17} for details) 
that there are $\delta_0 >0$ and open sets $\mathcal{M}_n$, $\mathcal{K}_n$ and $\mathcal{L}_n$, for $n \geq 0$, 
such that
\begin{itemize}
\item[(i)] $Y_n' \subset \mathcal{M}_n$, $(Y_n' \setminus  Y_{n,\diamond}) \subset \mathcal{K}_n$, and 
$Y_{n, \diamond} \subset \mathcal{J}_n$, $\mathcal{M}_n=\mathcal{K}_n \cup \mathcal{J}_n$, 
\item[(ii)] for all $j \in \J_n$, $B_{\delta_0}(\chi_{n,j}(\mathcal{M}_n)+j) \subset \mathcal{M}_{n-1}$, and 
$B_{\delta_0}(\chi_{n,a_n-1}(\mathcal{J}_n)) \subset \mathcal{M}_{n-1}$. 
\end{itemize}
Compare (ii) with Lemma~\ref{lemma:tiled}-(a). 
%Also, note that because $\Lambda_n= (\chi_{n+1}(\Lambda_{n+1}) + \Z) \cap Y_n'$, it 
It follows that $B_{\delta_0}(\Lambda_n) \subset \mathcal{M}_n$.
  
Let $\rho_n(\zeta)|d\zeta|$ denote the hyperbolic metric of constant curvature $-1$ on $\mathcal{M}_n$. 
By a proof similar to the one of the Schwrz-Pick lemma (see \cite[Lemma 5.5]{Che19} or \cite[Lemma 3.8]{AC18} 
for more details) one can show the following.

\begin{cor}\label{cor:hyperbolic-contraction}
There exists $\mu < 1$ such that for all $ n \geq 0 $, the following hold: 
\begin{itemize}
\item for all $ j \in \J_{n-1}$ and all $\zeta \in Y_n'$, $|\chi_{ n,j}' (\zeta) |\,\rho_{ n-1 } \p{ \chi_{ n,j} (\zeta) } \leq \mu\,\rho_{ n } (\zeta)$, 
\item for all $\zeta \in Y_{n, \diamond}'$, $|\chi_{ n,a_n-1}' (\zeta) |\,\rho_{ n-1 } \p{ \chi_{ n,a_n-1} (\zeta) } \leq \mu\,\rho_{ n } (\zeta)$. 
\end{itemize}
\end{cor}

For each $\zeta_0 \in \Lambda_0 \cup \Delta_0$, there are unique sequences 
\[\ea(\zeta_0)=(s_i)_{i\geq 1}, \qquad  O_r(\zeta_0):=(\zeta_i)_{i\geq 0},\]
such that for all $i\geq 0$, $s_i \in \Z\cap [0, a_i-1]$, $\zeta_i \in \Lambda_i \cup \Delta_i$, 
and $\chi_{i+1}(\zeta_{i+1}) + s_{i+1}= \zeta_i$.  
The sequence $\ea(\zeta_0)$ is called the \textbf{address} of $\zeta_0$, and $O_r(\zeta_0)$ is called 
the \textbf{renormalisation orbit} of $\zeta_0$. 
%In terms of the notations in the above paragraph, for each $n\geq 1$, 
%$\xi_n: \Lambda_{n-1} \cup \Delta_{n-1} \to \Lambda_n \cup \Delta_n$ is defined as $\xi_n(\zeta_{n-1})= \zeta_n$.
For an address $\ea=(s_n)_{n\geq 1}$, and integers $0 \leq m \leq n$, we define the notations 
\[\chi_{ n \to m,\,\textbf{s}} = \chi_{ m + 1 , s_{ m + 1 } } \circ \cdots \circ \chi_{ n, s_n },\]  % \xi_{m\to n}=\xi_n\circ\cdots\circ\xi_{m+1},\]
with the convention that when $ m = n $, $ \chi_{ n \to n,\,\ea}$ is the identity map.

\subsection{Renormalisation orbits of hairs}
In this section we identify a key feature of the renormalisation orbits of points on the set of hairs.

\begin{lemma} \label{lema:in-Siegel-disk}
For any $\zeta_0 \in \Lambda_0$ with renormalisation orbit $O_r(\zeta_0)= (\zeta_i)_{i\geq 0}$,
if there is $D_0>0$ such that $\im \zeta_n > D_0/ \alpha_n$ for all $n \geq 0$, then $\zeta_0 \in \overline{\Delta}_0$.
\end{lemma}

Lemma~\ref{lema:in-Siegel-disk} places a restriction on the renormalisation orbit of the points in the hairs. 
That is, if $ \zeta_0 \in \Lambda_0\setminus\overline{\Delta}_0$, then for any $D_0>0$ there are infinitely many levels $n$
at which $\zeta_n$ lies below the height $D_0/\alpha_n$. 

\begin{proof}
Fix arbitrary $\zeta_0$ and $D_0$ as in the hypothesis. 
We claim that there is a constant $M>0$ such that if $\zeta'_0 \in \Lambda_0\cup\Delta_0$ satisfies
$\im \zeta'_0 \geq \im \zeta_0 + M$, then $O_r (\zeta_0')=(\zeta'_n )_{n \geq 0}$ satisfies 
$\im \zeta'_n \geq \im \zeta_n+M $ for all $ n \geq 0 $. 
%  Define as well $ y'_n = \operatorname{Im} y_n $.

Because the width of $Y_n'$ is comparable to $1/\alpha_n$, there is a constant $D_0'>0$, independent of $n$, such that 
for any $\zeta \in Y_n'$ satisfying $\im \zeta \in [-2, D_0/\alpha_n]$, 
\begin{equation}\label{equ:est-4}
\log(1+|\zeta|)\leq \log (1/\alpha_n)+D_0' \text{\quad and\quad}
\log\left( 1+\left | \zeta-1/\alpha_n \right|\right) \leq \log (1/\alpha_n)+{D_0'}.
\end{equation}
Let $M_0=M_0(D_0)$ and $\widetilde{M}_0=\widetilde{M}_0(D_0)$ be the corresponding constants 
in Proposition~\ref{P:est-imag}. Let  
\begin{equation}\label{equ:M-range}
M:= \max \left \{\widetilde{M}_0+M_0 + D_0'/(2\pi), 4M_0 \right\}+1.
\end{equation}
Now we prove by induction that $\im \zeta_n' \geq \im\zeta_n + M$  for all $ n \geq 1$. 
Assume that for some $n \geq 1$, $\im \zeta_{n-1}' \geq \im \zeta_{n-1} + M$. 
If $\im \zeta_n' < D_0/\alpha_n$, by \eqref{equ:est-4} and Proposition~\ref{P:est-imag}, we must have
\[\im \zeta_{n-1}' \leq \frac{1 }{ 2 \pi } \log \frac{1}{\alpha_{n}} +\frac{D_0}{2\pi}+ \widetilde{M}_0, \qquad 
\im\zeta_{n-1} \geq \frac{1 }{ 2 \pi } \log \frac{1}{\alpha_n}+D_0 -M_0.\]
These contradict the induction hypothesis. %$\im \zeta_{n-1}' \geq \im \zeta_{n-1}+M$. 
Thus, $ \im \zeta_n' \geq  D_0/\alpha_n$, and by Proposition~\ref{P:est-imag}-(a),
\[\begin{split}
\im \zeta_n' \geq & \frac{1}{\alpha_n}\im\zeta_{n-1}'-\frac{1}{2\pi\alpha_n} \log \frac{1}{\alpha_n}-\frac{M_0}{\alpha_n}\\
\geq & \left(\frac{1}{\alpha_n}\im\zeta_{n-1}-\frac{1}{2\pi\alpha_n}\log\frac{1}{\alpha_n}
+\frac{M_0}{\alpha_n}\right)+\frac{M-2M_0}{\alpha_n}
\geq \im\zeta_n+M.
\end{split}\]

Consider the sets $\Omega_n = \{\zeta \in Y_n' \mid \im \zeta \geq \im \zeta_n + M\} $, for $n \geq 0$. 
By the above paragraphs, $\Omega_{n-1}$ projects into $\Omega_n$ via $\chi_{n, j}^{-1}$. 
Then, by the definition of renormalisation, it follows that each $f_n$ can be iterated infinitely many times on 
$\Phi_n^{-1} (\Omega_n)$. 
Therefore, $\Omega_n \subseteq \Delta_n$, for all $n\geq 0$. 
In particular, each $\zeta_n$ is uniformly close to some $u_n \in \Delta_n$. 
Thus, by Corollary~\ref{cor:hyperbolic-contraction}, $\zeta_0 \in \partial \Delta_0$. 
%By the definition of $\MD_n$, there exist a constant $C>0$ and a sequence of real numbers $(x_n)_{n\in\N}$ with $|x_n|\leq C$ such that $ u_n = \zeta_n + \ii M +x_n\in\MD_n$ for all $n\in\N$.  By following the same address as $ \zeta_n $ and pulling it upward to the level $ 0 $ of the renormalization tower, we obtain a point $ w_n \in \MD_0\cap\Delta_0 $ for each $n\in\N$. It follows from Corollary \ref{cor:converge-seq-pt} that $ w_n \to \zeta_0 $ as $n\to\infty$. Therefore we have $\zeta_0\in\overline{\Delta}_0$.
\end{proof}

\begin{lemma} \label{lema:evetually-below}
Assume that $\alpha \in  \HJ \cup \HS$. 
For any $\zeta_0 \in \Lambda_0 \setminus \overline{\Delta}_0$ with $O_r(\zeta_0)= (\zeta_n)_{n\geq 0}$ 
and any $D_0>0$, there exists $n_\star \geq 0 $ such that for all $n \geq n_\star$, 
\begin{equation} \label{equ:point-below-exp-model}
\im \zeta_n \leq D_0/\alpha_n. 
\end{equation}
\end{lemma}

\begin{proof}
Fix arbitrary $\zeta_0$ and $D_0>0$, and let $M_0=M_0(D_0)$ be the corresponding constant from Proposition~\ref{P:est-imag}.
We claim that, for each $\alpha \in  \HJ \cup \HS$, there exists $n_0 \geq 0$ such that for all $n \geq n_0$, 
\begin{equation} \label{equ:fast-growth-rot}
\frac{ 1 }{ 2 \pi } \log \frac{ 1 }{ \alpha_{ n + 1 } }>\frac{ D_0 }{ \alpha_n } + M_0 - D_0.
\end{equation}
By a direct calculation, and using $\log(1-x)\geq -2x$ for $0\leq x\leq 1/2$, for $\alpha \in  \HJ$, 
\[
\frac{ 1 }{ 2 \pi } \log \frac{ 1 }{ \alpha_{ n + 1 } }\geq \frac{ 1 }{ 2 \pi } \log \left(a_{n+1}-\frac{1}{2}\right)
\geq \frac{ 1 }{ 2 \pi } \log a_{n+1}-\frac{1}{2\pi}\frac{1}{a_{n+1}} 
> a_n\,\frac{ u_n \log a_n}{2\pi } -1\geq \frac{u_n\log a_n}{4\pi\alpha_n}-1.\]
Because $ u_n \log a_n \to +\infty $ as $n\to\infty$, there is $n_1'\geq 0$ such that for all 
$n \geq n_1'$ \eqref{equ:fast-growth-rot} holds.

Now assume that $\alpha \in \HS$ with the corresponding sequences $(\eta_n)_{n \geq 0}$ and $(v_n)_{n \geq 0}$. 
Suppose that $|\eta_n| \leq C'$ for all $n \geq 0$ and $e^{v_n} \geq 2C'$ for all $n \geq n_2$. 
Then, for all $n\geq n_2$, 
\[
\frac{ 1 }{ 2 \pi } \log \frac{ 1 }{ \alpha_{ n + 1 } } 
\geq \frac{ 1 }{ 2 \pi } \log a_{n+1}-\frac{1}{2\pi}\frac{1}{a_{n+1}} 
> \frac{ 1 }{ 2 \pi } \log \Big(e^{v_n a_n}+\eta_n\Big)-1 
\geq \frac{ v_n }{ 4 \pi\alpha_n } -\frac{C'}{\pi}-1.\]
Since $v_n\to +\infty$ as $n\to\infty$, there exists $n_2'\geq n_2$ such that for all $n\geq n_2'$, 
\eqref{equ:fast-growth-rot} holds. 
This completes the proof of the claim. 
%Therefore, if $n\geq n_0=\max\{n_1',n_2'\}$, for all $\alpha \in \HJ\cup\HS$, we have \eqref{equ:fast-growth-rot}.

By Lemma~\ref{lema:in-Siegel-disk}, there exists $ n_\star\geq n_0+1$ such that 
$\im \zeta_{n_\star}\leq D_0/\alpha_{n_\star}$. 
We claim that if $\zeta_n$ satisfies \eqref{equ:point-below-exp-model} for some $ n \geq n_\star $, 
then $\zeta_{ n + 1 }$ also satisfies \eqref{equ:point-below-exp-model}. 
Otherwise, if in the contrary $ \im\zeta_{ n + 1 } > { D_0 }/{ \alpha_{ n + 1 } } $, then 
by Proposition~\ref{P:est-imag}(a) and \eqref{equ:fast-growth-rot}, we get the contradiction 
\[\im\zeta_n \geq \alpha_{n+1} \im\zeta_{ n + 1 } + \frac{ 1 }{  2 \pi } \log \frac{ 1 }{ \alpha_{ n + 1 } } - M_0
>\frac{ 1 }{  2 \pi } \log \frac{ 1 }{ \alpha_{ n + 1 } }+D_0-M_0>\frac{D_0}{\alpha_n}.\qedhere\]
\end{proof}

For $\kappa > 1$, and $n\geq 0$, define 
\[\MD_n^\kappa = \left \{\zeta \in Y_n' \mid  
\im \zeta \geq \min \left\{\abs{ \re\zeta }^{\kappa}, \abs{\re \zeta - 1/\alpha_n}^{\kappa}\right\}\right\}.\]

%a point $ \zeta \in \MD_n $ is said to be \emph{above the $\kappa $-parabola (in $\MD_n $)} if it satisfies the following:
%  curves y = x^p are called ``higher parabolas''
%\[\im\zeta \geq \abs{ \re\zeta }^{\kappa}  \text{ or } \im\zeta \geq \abs{ \frac{ 1 }{ \alpha_n } - \re\zeta }^{\kappa}.\]
%The set of points above the $\kappa $-parabola in $\MD_n $ will be denoted by $$.
%Similarly we will say that a point $ \zeta \in \MD_n $ is \emph{below the parabola} if it is not above the parabola.

%For $n \geq 0$, a point $\zeta \in \Lambda_n $ is called \textbf{accessible} if it is an end point of an (open) arc in 
%$Y_n' \setminus \Lambda_n$. 
%By the definition of Cantor bouquet and one-sided hairy circle (see Section~\ref{subsec:topo}), and 
%Theorem~\ref{thm-Lambda-topo}, if $\alpha\in\HT_N\setminus\HH$, then each component of 
%$\Lambda_n \setminus \overline{\Delta}_n$ is accumulated from both sides by components of 
%$\Lambda_n\setminus\overline{\Delta}_n$.
%This means that the set of accessible points of $\Lambda_n \setminus \overline{\Delta}_n$ 
%is precisely the set of end points $\ME_n$. 

\begin{lemma} \label{lem:parabola-to-accessible}
Assume that $\alpha \in  \HJ \cup \HS$. If $\zeta_0 \in \Lambda_0\setminus\overline{\Delta}_0$ with 
$O_r(\zeta_0)=(\zeta_n)_{n\geq 0}$ satisfies at least one of the following properties, 
\begin{itemize}
\item[(i)] there are $\kappa >1 $ and infinitely many distinct $n$ such that $\zeta_n \notin \MD_n^\kappa$, 
\item[(ii)] $\liminf_{n \to \infty} \im \zeta_n < \infty$,
\end{itemize} 
then $\zeta_0 \in \ME_0$.
\end{lemma}

\begin{proof}
Let $v_0 \in \ME_0$ denote the end point of the component of $\Lambda_0\setminus \overline{\Delta}_0$ 
which contains $\zeta_0$. Using $\ea(\zeta_0)= (s_i)_{i \geq 1}$, we have 
$v_0 = \lim_{n\to\infty} \chi_{ 1, s_1 } \circ \cdots \circ \chi_{ n, s_n } (1) $. 
We aim to show that either of (i) and (ii) implies that $\zeta_0=v_0$. 

Let us first assume that (i) holds. 
There are $n_1 < n_2< n_3 < \dots$ and $\kappa >1$ such that for all $j\geq 1$, 
\begin{equation}\label{equ:zeta-below-para}
\im \zeta_{ n_j } < | \re \zeta_{ n_j } |^{\kappa} \text{\quad and\quad} 
\im \zeta_{ n_j } < | 1 / \alpha_{ n_j } - \re \zeta_{ n_j } |^{\kappa}.
\end{equation}
Without loss of generality, for convenience, we assume that $\re\zeta_{n_j}\leq 1/(2\alpha_{n_j})$ for all $j\geq 0$. 
The other cases are dealt with by symmetry. 

Applying Lemma~\ref{lema:evetually-below} with $D_0=1$, there is $n_\star \geq 0$ such that 
$\im \zeta_n \leq 1/\alpha_n$ for all $ n \geq n_\star$. 
For convenience, let us assume that $n_1 \geq n_\star + 2$. 
%  Since $ \alpha $ is not Hermann $ \zeta_0 $,
%  there is $ n $ such that $ \operatorname{Im} \zeta_n \leq \frac{ 1 }{ 2 \pi } \log \frac{ 1 }{ \alpha_n } + D_4 + 2 $.
%  Without loss of generality we assume that
%  $ \operatorname{Im} \zeta_0 \leq \frac{ 1 }{ 2 \pi } \log \frac{ 1 }{ \alpha_n } + D_4 + 2 $.

Recall that $\gamma_n \cup \tilde{\gamma}_n$ is contained in $\mho$, and hence all of its integer translate are pairwise 
disjoint.  
Let us choose a sequence of real numbers $(x_n)_{n \geq 0}$ in $(0,1)$ such that for all $n \geq 0$, 
$\re\zeta_n +x_n \in Y_n' \cap \Z$, and define 
%\[\re\zeta_n +x_n \in 
%\begin{cases}
%Y_n' \cap \Z \cap [1,1/(2\alpha_n)]   & \text{ if }   \re\zeta_n\leq 1/(2\alpha_n),\\
%Y_n' \cap \Z \cap (1/(2\alpha_n),+\infty)   & \text{ if }  \re\zeta_n> 1/(2\alpha_n).
%\end{cases}.\]
\begin{equation}\label{equ:u-j}
u_n =\re \zeta_n + x_n \text{\quad and\quad} u_{n-1}'= \chi_{n, s_n}(u_n).
\end{equation}

Because $\im\zeta_{ n_j }$ and $\im u_{n_j}$ are at most $1/\alpha_{n_j}$, applying Proposition~\ref{P:est-imag}-(b) 
with $D_0=1$, 
\[\im\zeta_{n_j-1} \geq \frac{1}{2\pi}\log \left( 1 + | \zeta_{ n_j } | \right)-\widetilde{M}_0,  \qquad 
\im u_{n_j-1}' \geq \frac{1}{2\pi}\log \left( 1 +  u_{n_j}  \right)-\widetilde{M}_0,\]
and applying Proposition~\ref{prop:esti-chi}-(b) with $D_0=1$, 
\begin{align*}
|\zeta_{n_j-1} - u_{n_j-1}'| 
& \leq \int_{u_{n_j}+(\zeta_{n_j}- u_{n_j}) [0,1]}  |\chi_{n_j}'(z)| \, |dz| \\
%& \leq \int_0^{\zeta_{n_j} - u_{n_j}} \frac{1}{|u_{n_j}+ t|} \, |dt|\leq \log |\zeta_{n_j}| - \log |u_{n_j}| \\
& \leq \widetilde{M}_1 \log |\zeta_{n_j} - u_{n_j}| 
\leq \widetilde{M}_1 \log |u_{n_j}|^\kappa= \widetilde{M}_1 \kappa \log |u_{n_j}|.  
\end{align*}
Similarly, repeating the integral estimate, and using the above two equations, we conclude that 
\begin{equation}\label{E:bounded-images}
\left| \zeta_{n_j-2}- \chi_{ n_j - 1, s_{ n_j - 1 } } (u_{n_j-1}')\right| = O(1)
\end{equation}
with the constant in $O$ independent of $j$. 
On the other hand, $|\chi_{ n_j+1, s_{ n_j+1 } } (1)- u_{n_j} |$ is also uniformly bounded from above. 
Then, by Corollary~\ref{cor:hyperbolic-contraction}, %we conclude that 
\begin{align*}
\zeta_0= \lim_{j \to \infty} \chi_{ 1, s_1 } \circ \cdots \circ \chi_{n_j, s_{ n_j } } (\zeta_{n_j}) 
& = \lim_{j \to \infty} \chi_{ 1, s_1 } \circ \cdots \circ \chi_{n_j, s_{ n_j } } (u_{n_j}) \\
& = \lim_{j \to \infty} \chi_{ 1, s_1 } \circ \cdots \circ \chi_{ n_j, s_{ n_j } }\circ \chi_{ n_j+1, s_{n_j+1}}(1) = v_0.
\end{align*}

If (ii) occurs, $|\zeta_{n} - u_n|$ is uniformly bounded from above for infinitely many $n$, 
and the final step of the argument applies.
\end{proof}

\subsection{Bounding the dimension of the hairs}
By making $\delta_0>0$ smaller if necessary, we may assume that $\delta_0 \in (0,1)$. 
For $n\geq 0$, consider the finite collection 
\[\MQ_n = \left\{Q_n = \Boxx\left(\frac{u\delta_0 +i v\delta_0}{4},\frac{\delta_0}{8}\right) \,\middle\vert\, 
u,v \in \N, Q_n \cap (\Lambda_n\cup\Delta_n) \neq \emptyset \right\}.\]
For each $\zeta \in \Lambda_n \cup \Delta_n$ with $\im \zeta>0$, there is $Q_n(\zeta) \in \MQ_n$ that contains $\zeta$. 
%Also, because the diameter of each $Q_n \in \MQ_n$ is less than $\delta_0/2$, by the choice of $\delta_0$, each 
%$\chi_n$ is defined on every $Q_n \in \MQ_n$. 
%Although, $\chi_n$ may not be univalent on some $Q_n \in \MQ_n$, if $Q_n$ has sufficiently large height,
%$\chi_n$ has univalent extension onto $B_{\delta_0/2}(Q_n)$. 
For $n\geq 0$, let 
\[\MK_n=
\left\{ \chi_{n\to 0, \ea}(Q_n) \,\middle\vert\, 
\forall n\geq 1, Q_n \in \MQ_n, \ea=(s_1,s_2, s_3,\cdots) \text{ with } s_n \in \widetilde{\J}_{n-1} 
\right\}.\]

\begin{lemma}\label{lema:exp-behavior}
There is a constant $C > 0$ such that for any $\kappa >1$ and any $n\geq 0$ the following hold:
\begin{enumerate}
\item for any $Q_n \in \MQ_n$ which lies in the region $\im \zeta \geq C$, $\chi_n$ has univalent extension onto 
$B_1(Q_n)$, 
\item for any $\zeta \in B_1(\MD_n^\kappa)$ with $C \leq \im \zeta \leq 2/\alpha_n$, 
$\im \zeta \geq e^{\im \chi_n(\zeta)}$. 
\end{enumerate}
\end{lemma}

\begin{proof}
(a) The map $\chi_n$ is univalent on $Y_n'$, and has a unique critical point, say $c \in \partial Y_n'$. 
By the definition of renormalisation, $\Expo(c)$ is the unique critical point of $\MMR (f_n)$. 
On the other hand, $\MMR f_n=P \circ \psi^{-1}$, for some univalent map $\psi^{-1}: U \to \C$ 
with $\psi(0)=0$ and $|\psi'(0)|=1$. 
It follows from the Koebe distortion theorem that the critical point of $\MMR f_n$, $\psi^{-1}(-1/3)$, 
is uniformly away from $0$. 
Therefore, $\im c \leq C'$, for a uniform $C'$. 
Thus, if $Q_n$ lies above the line $\im \zeta = C + 1$, $\chi_n$ has a univalent extension onto $B_1(Q_n)$. 

(b) Fix an arbitrary $\zeta_n \in B_{1}(\MD_n^\kappa)$, and let $\zeta_{n-1}= \chi_n(\zeta_n)$. 
Note that if $\im \zeta_n \geq 1$, 
\begin{equation}\label{E:lem:expo-behavior}
(1/2) \min\{|\zeta|,|\zeta-1/\alpha_n|\} \leq \im\zeta \leq \min\{|\zeta|,|\zeta-1/\alpha_n|\}.
\end{equation}
By Proposition~\ref{P:est-imag}-(b), there is $\widetilde{M}_0$ independent of $n$ such that  
%\[\big|\im \zeta_{n-1}-\tfrac{1}{2\pi}\log(1+|\zeta_n|)\big|\leq \widetilde{M}_0.\]
%If $\im\zeta_{n-1}\geq\widetilde{M}_0+1$, then we have
\[2\pi(\im\zeta_{n-1}-\widetilde{M}_0)\leq\log(1+|\zeta_n|)\leq 2\pi(\im\zeta_{n-1}+\widetilde{M}_0).\]
By symmetry, we may assume that $\re \zeta_n \leq 1/(2\alpha_n)$. 
Thus, combining the above equations, we get 
\begin{equation}\label{E:lema:exp-behavior}
(1/2)(e^{2\pi \im\zeta_{n-1}} e^{-2\pi \widetilde{M}_0} -1) 
\leq \im \zeta_n \leq e^{2\pi \im\zeta_{n-1}} e^{+2\pi \widetilde{M}_0}-1.
\end{equation}
There is a constant $C'$, independent of $n$, such that if $\im \zeta_{n-1} \geq C'$, the left hand side of the above 
equation is at least $e^{\im \zeta_{n-1}}$. 
The right hand inequality in the above equation implies that if $\im \zeta_n$ is sufficiently large, 
$\im \zeta_{n-1} \geq C'$. 
\end{proof}

With the constant $C$ from Lemma~\ref{lema:exp-behavior}, $\kappa>1$ and $n\geq 0$, let  
\[\MQ_n^\kappa= \{Q_n \in \MQ_n \mid Q_n \cap \MD_n^\kappa \neq \emptyset, \text{ and } 
\forall \zeta \in Q_n, C \leq \im \zeta \leq 2/\alpha_n\}.\]

%For $\kappa >1$, and $C>0$, we define
%\[\MD_n^\kappa(\delta_0, C) = \{\zeta \in B_{\delta_0/2}(\MD_n^\kappa)  \mid C \leq \im\zeta \leq 1/\alpha_n+1\}.\]

%\[\Pi_n^\kappa(D_0,C) = B_{\delta_0/2}(\Lambda_n\cup\Delta_n) \bigcap 
%\{\zeta \in  \MD_n^\kappa \mid C\leq\im\zeta\leq D_0/\alpha_n+1\}.\]

By the definition of $\HJ$ and $\HS$, for large enough $n$,  $C < 1/\alpha_n$, and hence 
$\MQ_n^\kappa \neq \emptyset$. 

For $Q_{n-1}\in \MQ_{n-1}$, let 
\[I(Q_{n-1}) = \max \{\im \zeta \mid \zeta \in Y_n', \exists i \in \N, \chi_{n,i}(\zeta) \in Q_{n-1}\}.\]

\begin{lemma}\label{lema:exp-behavior-2}
There is a constant $M_2$ such that for any $\kappa >1$ and any $n\geq 0$ the following holds. 
If $Q_{n+1} \in \MQ_{n+1}^\kappa$ and $Q_{n} \in \MQ_{n}$ satisfy $\chi_{n+1, i}(Q_{n+1}) \cap Q_{n} \neq \emptyset$ for some $i \in \N$, 
then for all $\zeta\in Q_{n+1}$, 
\begin{equation}\label{E:lem:expo-behavior-2}
\frac{1}{M_2 I(Q_{n})}\leq |\chi_{n+1,i}'(\zeta)|\leq \frac{M_2}{ I(Q_{n})}.
\end{equation}
\end{lemma}

\begin{proof}
Fix an arbitrary $\zeta_{n+1} \in Q_{n+1}$ with $\zeta_n= \chi_{n+1,i}(\zeta_{n+1}) \in Q_n$, 
and let $\zeta_{n+1}' \in Y_{n+1}'$ realises the $\max$ in the definition of $I(Q_n)$.  
In particular, $\im \zeta_{n+1}' \geq \im \zeta_{n+1}$. 
By Proposition~\ref{prop:esti-chi}(b), there is $\widetilde{M}_1\geq 1$, 
such that 
\[\widetilde{M}_1^{-1}/|\zeta_{n+1}| \leq |\chi_{n+1,i}'(\zeta_{n+1})| \leq \widetilde{M}_1/|\zeta_{n+1}|.\]
By \eqref{E:lem:expo-behavior}, $|\zeta_{n+1}|/\im \zeta_{n+1}$ is bounded from above and below. 
Thus, the lower bound in \eqref{E:lem:expo-behavior-2} holds.

Because $|\im \zeta_n - \im \zeta_n'| \leq \delta_0/4$, $\im \zeta_{n+1}'= O (\im \zeta_{n+1})$, and 
also by the Koebe distortion theorem, $|\chi_{n+1, i}'(\zeta_{n+1}')|/|\chi_{n+1, i}'(\zeta_{n+1})|$ is uniformly bounded 
from above and below. 
\end{proof}

%Recall that 
%\[H_{n-1}= 
%\big(\bigcup_{ j \in \mathbb{J}_{n-1}} \chi_{ n, j } (H_n)\Big)\cup \chi_{ n, J_n } (H_n\setminus H_{n,\diamond})
%\subset \bigcup_{ j \in \widetilde{\J}_{n-1} } \chi_{ n, j } ( H_n ).\]

% \begin{lemma} \label{lemma: above parabola have dimension 1}
%  SOUNDS LIKE THIS IS GOING TO BE THE MAIN STATEMENT
%  HERE...
% \end{lemma}

Recall the classes $\HJ$ and $\HS$ from Section~\ref{section:arithmetic rot numbers}, 
and the class of maps $\QIS_\alpha$ from Section~\ref{subsec-IS-class}. 

\begin{thm}\label{thm:dim-hairs}
For every $\alpha \in \HJ \cup \HS$ and every $f \in \QIS_\alpha$, we have
\[\dim_H \big (\Lambda(f) \setminus (\overline{\Delta}(f) \cup \ME(f)) \big)=1 
\text{\quad and\quad}\dim_H(\ME(f))=2.\]
\end{thm}

\begin{proof}
It is enough to show that for any $\varepsilon \in (0,1)$ and any $Q_0 \in \MQ_0$, $\dim_H (H_0\cap Q_0) \leq 1+\varepsilon$. 
Fix an arbitrary $\varepsilon \in (0,1)$ and $Q_0 \in \MQ_0$. Define $\kappa =2/\varepsilon$. 

For any $n \geq 0$ and $\zeta \in \Lambda_n$, let $Q_n (\zeta)$ be an element of $\MQ_n$ which contains $\zeta$. 
That may not be unique, but any such choice works in the sequel. 

For $\zeta_0 \in \Lambda_0$ consider the renormalisation orbit $O_r (\zeta_0)= (\zeta_i)_{i \geq 0}$ down the tower, 
and the address $\ea= \ea(\zeta_0)$. 
Let $C$ be the constant from Lemma \ref{lema:exp-behavior}. 
For $k \geq 0$, let $V_k$ denote the set of all $\zeta_0 \in \MH_0 \cap Q_0$ such that for all $n \geq k$, 
\[ \zeta_n \in \MD_n^\kappa, \qquad \text{ and } \qquad C+1 \leq \im \zeta_n \leq 1/\alpha_n.\]
By Lemmas \ref{lema:evetually-below} and \ref{lem:parabola-to-accessible}, $Q_0 \cap H_0  \subset \cup_{ k \geq 0}  V_k$.
Thus, it suffices to show $\dim_H (V_k) \leq 1 + \varepsilon $ for all $ k \geq 0$. 

Let us first assume that $k = 0$. 
For $n \geq k$, let $\MA_n$ be the collection of all $\chi_{n \to 0, \ea(\zeta_0)}(Q_n(\zeta_n))$ over all $\zeta_0\in V_k$. 
Each $\MA_n$ is a covering of $V_k$, and by Corollary~\ref{cor:hyperbolic-contraction}, 
\[\lim_{n \to \infty} \max_{K_n\in\MA_n}\diam\,K_n=0.\]
It is enough to show that  
\begin{equation}\label{equ:Hdim-measure}
\sup_{n\geq 0} \textstyle{\sum_{K_n\in\MA_n}} (\diam\, K_n)^{1+\varepsilon} < \infty.
\end{equation}
For $K_n\in\MA_n$, let $\MG(K_n)$ denote the collection of all $K_{n+1} \in \MA_{n+1}$ such that 
$K_{n+1} \cap K_n \neq \emptyset$. 
In order to prove \eqref{equ:Hdim-measure}, it suffices to prove that there exists $n_0\geq k$ such that for all $n\geq n_0$ 
and all $K_n\in\MA_n$,
\begin{equation}\label{equ:criterion}
\textstyle{\sum_{K_{n+1} \in \MG(K_n)}} (\diam\, K_{n+1})^{1+\varepsilon}\leq (\diam\,K_n)^{1+\varepsilon}.
\end{equation}

Let $n\geq k$ and fix an arbitrary $K_n=K_n(\zeta_0) \in \MA_n$.
For $0 \leq i \leq n$, let $Q_i= Q_i(\zeta_i)$, and define $I_{i+1}=I(Q_i)$. 

By Lemma~\ref{lema:exp-behavior}-(a) and the Koebe distortion theorem, there exists a constant $C_1 \geq 1$ 
independent of $n$ such that for any $\zeta \in Q_n$, we have
\begin{equation}\label{equ:distortion-bd}
C_1^{-1}|\chi_{n\to 0}'(\zeta)|\leq \diam\,K_n(\zeta_0) \leq C_1|\chi_{n\to 0}'(\zeta)|,
\end{equation}
and hence by Lemma \ref{lema:exp-behavior-2}, we have
\[\textstyle{ \diam K_n(\zeta_0) \geq\frac{1}{C_1}|\chi_{n\to 0}'(\zeta)|
=\frac{1}{C_1} \prod_{i=1}^n |\chi_i'(\chi_{n\to i}(\zeta))|
\geq \frac{1}{C_1}\prod_{i=1}^n \frac{1}{M_2 I_i}}.\]
On the other hand, any $K_{n+1}\in\MG(K_n)$ can be written as 
$K_{n+1}=\chi_{n+1\to 0, \ea}(Q_{n+1})\in\MA_{n+1}$, where $Q_{n+1} \in \MQ_{n+1}^\kappa$ and 
$\chi_{n+1}(Q_{n+1}) \cap Q_n(\zeta_n) \neq\emptyset$.
For any $\zeta' \in Q_{n+1}$ (in particular when $\chi_{n+1}(\zeta') \in Q_n(\zeta_n)$), 
still by Lemma \ref{lema:exp-behavior}(b) and \eqref{equ:distortion-bd} we have
\[\textstyle{\diam K_{n+1}\leq C_1|\chi_{n+1\to 0}'(\zeta')|
=C_1\prod_{i=1}^{n+1} |\chi_i'(\chi_{n+1\to i}(\zeta'))|\leq C_1\prod_{i=1}^{n+1} M_2/I_i.}\]
The elements of $\MG(K_n)$ correspond to elements $Q_{n+1} \in \MQ_{n+1}^\kappa$ such that 
$\chi_{n+1}(Q_{n+1}) \cap Q_n \neq \emptyset$. Thus, the cardinality of $\MG(K_n)$ is at most 
$2 (4/\delta_0)^2 I_{n+1} I_{n+1}^{1/\kappa}$. 
Therefore, 
\[\textstyle{\sum_{K_{n+1}\in\MG(K_n)}\big(\diam\, K_{n+1}\big)^{1+\varepsilon}
\leq 2 (4/\delta_0)^2 I_{n+1}^{1+\varepsilon/2}\left(C_1\prod_{i=1}^{n+1}\frac{M_2}{I_i}\right)^{1+\varepsilon}
\leq \frac{\widetilde{C}_n}{I_{n+1}^{\varepsilon/2}}\big(\diam\,K_n\big)^{1+\varepsilon}.}\]
where $\widetilde{C}_n=2 (4/\delta_0)^2 C_1^{2(1+\varepsilon)}\,M_2^{(2n+1)(1+\varepsilon)}$.
By Lemma \ref{lema:exp-behavior}(a), $I_{n+1}$ increases exponentially fast, and eventually 
$I_{n+1}^{\varepsilon/2}\geq \widetilde{C}_n$. 
Thus, for large enough $n$, \eqref{equ:criterion} holds. 
We conclude that $\dim_H(V_0 \cap Q_0)\leq 1+ \varepsilon$. 

To deal with the other $V_k$ with $k \geq 1$, we note that there are a finite number of elements $Q_k$ in 
$\MQ_k$, and associated addresses $\ea$, such that $\chi_{k \to 0, \ea}(Q_k)$ covers $V_k$. 
By the above argument, each $Q_k \cap \MH_k$ has Hausdorff dimension at most $1+\varepsilon$. 
Since the changes of coordinates $\chi_{k \to 0, \ea}$ are holomorphic, we conclude that 
$\dim_h(V_k) \leq 1+ \varepsilon$.  
 
By Theorem~\ref{thm:dim-pc}, $\dim_H(\Lambda_0\setminus \overline{\Delta}_0 )=2$, and therefore, 
$\dim_H(\ME_0)=2$. 
\end{proof} 

For other families of maps enjoying Karpińska's dimension paradox see \cite{Sch07,Ber10}.
In the PhD thesis \cite{Bozovic2026}, the dimension paradox is established among infinitely satellite 
renormalisable quadratic polynomials where the post-critical set is a collection of arcs landing on a Cantor set 
of points (called hairy Cantor set in \cite{ChePe22}). 
%, and the dimension jumps from 0 to $r$ and then from $r$ to 2, 
%for the base Cantor set, hairs without end points, and all end points, respectively, 
%where $r \in [0,1]$ depends on the arithmetic of the satellite renormalisations. 

\subsection*{Acknowledgements.} 
D.~C.~acknowledges funding from EPSRC grant EP/M01746X/1. 
F.~Y. would like to thank Gaofei Zhang for providing financial support for a two months visit to Imperial College London
in 2018 (NSFC grant 11325104), and NSFC grant 12571093. 
%and the Fundamental Research Funds for the Central Universities in 2019 (grant 0203-14380025).

\bibliographystyle{amsalpha}
\bibliography{Data}
\end{document}